\documentclass[10pt]{article}

\newcommand{\sect}[1]{\section{#1}\setcounter{equation}{0}}

\title{Orthogonal Polynomials}
\def\shorttitle{Orthogonal Polynomials}

\author{Vilmos Totik}
\def\shortauthor{V. Totik}

\def\versiondate{11/11/05}

\def\abstracttext{In this survey, different aspects of
the theory of orthogonal polynomials of one
(real or complex) variable are reviewed.
Orthogonal polynomials on the unit circle are not discussed.}

\def\sep{\;\vrule\;}
\def\proof#1. {\par
                      \ifdim\lastskip<15pt
                      \removelastskip\penalty-200
                      \vskip15pt plus3pt minus3pt
                      \fi
                       {\def\a{#1}
                       \ifx\a\empty
                       {\noindent\bf Proof.}
                       \else
                       {\noindent\bf Proof of #1.}
                       \fi}\enspace}
\def\restr#1{\,\vrule\,\lower1.75ex\hbox{$#1$}}

\def\be{\begin{equation}}
\def\ee{\end{equation}}
\def\bea{\begin{eqnarray}}
\def\eea{\end{eqnarray}}
\def\bean{\begin{eqnarray*}}
\def\eean{\end{eqnarray*}}

\def\a{\alpha}
\def\b{\beta}
\def\d{\delta}

\def\e{\varepsilon}

\def\F{\Phi}
\def\g{\gamma}
\def\G{\Gamma}
\def\i{\infty}
\def\k{\kappa}
\def\l{\lambda}

\def\o{\omega}
\def\O{\Omega}

\def\r{\rho}
\def\s{\sigma}

\def\t{\theta}

\def\SZ{Szeg\H o}

\def\setm{\setminus}

\def\c{{\rm cap}}

\font\tenopen = cmbx10
\font\sevenopen = cmbx7
\font\fiveopen = cmbx5
\newfam\openfam
\def\open{\fam\openfam\tenopen}
\textfont\openfam = \tenopen
\scriptfont\openfam = \sevenopen
\scriptscriptfont\openfam = \fiveopen

\def\R{{\open R}}
\def\N{{\open N}}
\def\C{{\open C}}

\def\D{{\open D}}

\def\({\left(}
\def\){\right)}

\def\startpagenumber{70}
\def\volumenumber{1}
\def\year{2005}

\newcommand{\beginddoc}{
\begin{document}
\maketitle
\begin{abstract}
\abstracttext
\end{abstract}
\insert\footins{\scriptsize
\medskip
\baselineskip 8pt
\leftline{Surveys in Approximation Theory}
\leftline{Volume \volumenumber, \year.
pp.~\thepage--\pageref{endpage}.}
\leftline{Copyright \copyright\ 2005 Surveys in Approximation Theory.}
\leftline{ISSN 1555-578X}
\leftline{All rights of reproduction in any form reserved.}
\smallskip
\par\allowbreak}
\tableofcontents}
\renewcommand\rightmark{\ifodd\thepage{\hfill\it \shorttitle \hfill}\else {\hfill\it%
\shortauthor\hfill}\fi}
\renewcommand\leftmark{\ifodd\thepage{\hfill\it \shorttitle \hfill}\else {\hfill\it%
\shortauthor\hfill}\fi}
\markboth{{\it \shortauthor}}{{\it \shorttitle}}
\markright{{\it \shorttitle}}
\def\endddoc{\label{endpage}\end{document}}
\date{{\small \versiondate}}
\def\dword#1{{\bf #1}}
\beginddoc



\sect{Introduction}\label{secintro}

The theory of orthogonal polynomials can be divided into two main
but only loosely related parts. The two parts have many things in
common, and the division line is quite blurred, it is more or less
along algebra vs. analysis. One of the parts is the  algebraic
aspect of the theory, which has close connections with special
functions, combinatorics and algebra, and it is mainly devoted to
concrete orthogonal systems or hierarchies of systems such as the
Jacobi, Hahn, Askey-Wilson, $\ldots$ polynomials. All the discrete
polynomials and the $q$-analogues of classical ones belong to this
theory. We will not treat this part; the interested reader can
consult the three recent excellent monographs \cite{Ismail} by 
M. E. H. Ismail,
\cite{gautschi} by W. Gautschi and \cite{Askey} by 
G. E. Andrews, R. Askey and R. Roy. 
Much of the present state 
of the theory
of orthogonal polynomials of several variables lies also close to
this algebraic part of the theory. To discuss them would take us
too far from our main direction; rather we refer the reader to the
recent book \cite{xu} by C. F. Dunkl and Y. Xu.

The other part is the analytical aspect of the theory. Its methods
are analytical, and it deals with questions that are typical in
analysis, or questions that have emerged in and related to other
parts of mathematical analysis. General properties fill a smaller
part of the analytic theory, and the greater part falls into two
main and extremely rich branches: orthogonal polynomials on the
real line and on the circle. The richness is due to some special
features of the real line and the circle. Classical real
orthogonal polynomials, sometimes in other forms like continued
fractions, can be traced back to the 18th century, but their rapid
development occurred in the 19th and early 20th century.
Orthogonal polynomials on the unit circle are much younger, and
their existence is largely due to Szeg\H{o} and Geronimus in the
first half of the 20th century. B. Simon's recent treatise
\cite{simon1, simon2} summarizes and
greatly extends what has happened since then.

The connection of orthogonal polynomials with other
branches of mathematics is truly impressive. Without even trying to
be complete, we mention continued fractions, operator theory
(Jacobi operators), moment problems,
analytic functions (Bieberbach's conjecture),
interpolation, Pad\'e approximation, quadrature, approximation
theory, numerical analysis,
electrostatics, statistical quantum mechanics, special functions,
number theory (irrationality and transcendence),
graph theory (matching numbers), combinatorics,  random matrices,
stochastic processes (birth and death processes; prediction theory),
data sorting and compression, Radon transform and computer tomography.

This work is a survey on orthogonal polynomials
that do not lie on the unit circle. 
Orthogonal polynomials on the unit circle---both the classical theory
and recent contributions---will be hopefully dealt with in a companion article.

This work is meant for non-experts, and it therefore  contains
 introductory materials.
We have tried to list most of the actively researched fields not directly
connected with orthogonal polynomials on the unit circle, but because
of space limitation  we have only one or two pages on areas where dozens of papers
and several books had been published.  As a result,  our account is
necessarily incomplete. Also, the author's personal
taste and interest is reflected in the survey, and the omission of
a particular direction or a set of results reflects in no way
on the importance or quality of the omitted works.

For further backgound on orthogonal polynomials, the reader
can consult the books 
 Szeg\H{o} \cite{Szego}, Simon \cite{simon1}-\cite{simon2},
Freud \cite{freud},
Geronimus \cite{geronimus2}, Gautschi \cite{gautschi}, 
Chicara \cite{chihara}, Ismail \cite{Ismail}.
\medskip

This is a largely extended version of the first part of the article
Golinskii--Totik \cite{goltotik}.

\medskip

\noindent{\bf Acknowledgement.} Research was supported by  NSF grant
DMS-040650 and  OTKA  T049448, TS44782, and was carried out within
the Analysis Research Group of the Hungarian Academy of Sciences.
 The author participates in the
project INTAS 03-51-6637 that supported a visit by Leonid Golinskii
to  Szeged, during which the outline of this paper was laid
out.

\sect{Orthogonal polynomials}\label{secop}
\subsubsection{Orthogonal polynomials with respect to measures}
Let $\mu$ be a positive Borel measure on the complex plane, with 
an infinite number of points in its support,
for which 
\[\int |z|^md\mu(z)<\i\]
for all $m>0$. There are unique polynomials
\[p_n(z)=p_n(\mu,z)=\k_nz^n+\cdots,\qquad \k_n>0,\ n=0,1,\ldots\]
that form an orthonormal system in $L^2(\mu)$, i.e.
\[\int p_m \overline{p_n}d\mu=\left\{\begin{array}{ll}
  0 & \mbox{if $m\not= n$} \\
  1 & \mbox{if $m=n$.}
\end{array}\right.\]
These $p_n$'s are called the \dword{orthonormal polynomials}  corresponding to $\mu$.
$\k_n$ is the \dword{leading coefficient}, and $p_n(z)/\k_n=z^n+\cdots$
is called the \dword{monic orthogonal polynomial}. The leading coefficients
play a special and important role in the theory, many properties
depend on their behavior. When $d\mu(x)=w(x)dx$ on some interval, say,
then we talk about orthogonal polynomials with respect to the weight
function $w$.

The $p_n$'s can be easily generated:
all we have to do is to make sure  that
\be
\int \frac{p_n(z)}{\k_n}\overline{z^k}d\mu(z)=0,\qquad k=0,1,\ldots,n-1,
\label{system1}\ee
which is an $n\times n$ system of equations
for the non-leading 
coefficients of $p_n(z)/\k_n$ with matrix $(\s_{i,j})_{i,j=0}^{n-1}$,
where
\[\s_{i,j}=\int z^i\overline {z^j}d\mu(z)\]
are the \dword{complex moments} of $\mu$. This matrix is nonsingular:
if some linear combination with coefficients $c_0,\ldots,c_{n-1}$
of the rows is zero, then the polynomial $P_{n-1}(z)=c_0+\cdots+c_{n-1}z^{n-1}$
is orthogonal to every $z^j$, $j<n$, and hence it is orthogonal
to itself, i.e., 
\[\int|P_{n-1}|^2d\mu=\int P_{n-1}\overline{P_{n-1}}d\mu=0,\]
which implies $P_{n-1}(z)\equiv 0$. Thus, $c_0=\cdots=c_{n-1}=0$,
which shows the nonsingularity of $(\s_{i,j})$.
Therefore,  the system (\ref{system1}) 
has a unique solution for the non-leading coefficients of $p_n(z)/\k_n$
(note that the leading coefficient is 1), 
and finally $\k_n$ comes from normalization.

In particular, the complex moments already determine the polynomials.
In terms of them one can write up explicit determinant formulae:
\be p_n(z)=\frac{1}{\sqrt{D_{n-1}D_n}}\left|\begin{array}{ccccc}
  \s_{0,0} & \s_{0,1} & \cdots & \s_{0,n-1} & 1 \\
    \s_{1,0} & \s_{1,1} & \cdots & \s_{1,n-1} & z \\
  \vdots & \vdots & \ddots & \vdots & \vdots  \\
    \s_{n-1,0} & \s_{n-1,1} & \cdots & \s_{n-1,n-1} & z^{n-1} \\[5pt]
     \s_{n,0} & \s_{n,1} & \cdots & \s_{n,n-1} & z^n
\end{array}\right|\label{det}\ee
where
\be D_n=|\s_{i,j}|_{i,j=0}^n\label{gram}\ee
are the so called \dword{Gram determinants}.

Note that if $\mu$ is supported on the real line then
\[\s_{i,j}=\int x^{i+j}d\mu(x)=:\a_{i+j},\]
so $D_n=|\a_{i+j}|_{i,j=0}^n$ is a Hankel determinant, while if
$\mu$ is supported on the unit circle then
\[\s_{i,j}=\int z^{i-j}d\mu(z)=:\b_{i-j},\]
so $D_n=|\b_{i-j}|_{i,j=0}^n$ is a Toeplitz determinant. In these
two important cases the orthogonal polynomials have many special
properties that are missing in the general theory. 

\subsubsection{The Riemann--Hilbert approach}

Let $\mu$ be supported on the real line, and suppose that it is of the
form $d\mu(t)=w(t)dt$ with some smooth function $w$. A new approach to generating
orthogonal polynomials that has turned out to be of great importance was given in the
early 1990's by Fokas, Its and Kitaev \cite{Its}. Consider $2\times 2$ matrices
\[Y(z)=\left(\begin{array}{cc}
  Y_{11}(z) & Y_{12}(z) \\
  Y_{21}(z) & Y_{22}(z)
\end{array}\right)\]
where the $Y_{ij}$ are analytic functions on $\C\setm \R$, and solve
for such matrices the following matrix-valued \dword{Riemann--Hilbert problem}:

1. for all $x\in \R$
\[Y^+(x)=Y^-(x)\left(\begin{array}{cc}
  1 & w(x) \\
  0 & 1
\end{array}\right)\]
where $Y^+$, resp. $Y^-$, is the limit of $Y(z)$ as $z$ tends to $x$ from the upper, resp.
lower half plane, and

2. \[Y(z)=\left(I+O\left(\frac{1}{z}\right)\right)
\left(\begin{array}{cc}
  z^n & 0 \\
  0 & z^{-n}
\end{array}\right)\]
at infinity, where $I$ denotes the identity matrix.

Thus, one is looking for 4 functions $Y_{11},...,Y_{22}$
analytic on $\C\setm \R$, such that if $Y_{ij}^\pm(x)$
denote the boundary limits of these functions at $x\in \R$ from
the upper, resp. lower half plane, then
\be Y_{11}^-(x)=Y_{11}^+(x), \quad Y_{21}^-(x)=Y_{21}^+(x)\label{hil1}\ee
and
\be Y_{12}^-(x)=Y_{11}^+(x)w(x)+Y_{12}^+(x), \quad Y_{22}^-(x)=Y_{21}^+(x)w(x)+Y_{22}^+(x)\label{hil2}.\ee
These connect the functions on the upper and lower
half planes only very mildly, and what puts the problem
into a rigid form is the second condition, namely it is required
that for large $z$ uniformly on the plane we have
\be Y_{11}(z)=z^n+O(|z|^{n-1}), \quad  Y_{21}(z)=O(|z|^{n-1})\label{hil3}\ee
and
\be Y_{12}(z)=O(z^{-n-1}), \quad  Y_{22}(z)=z^{-n}+O(|z|^{-n-1}).\label{hil4}\ee
It can be shown that there is a unique solution $Y(z)$.
The relevance of this to orthogonal polynomials is
that the entry $Y_{11}(z)$
is precisely the monic polynomial $p_n(\mu,z)/\k_n$.  
Indeed, (\ref{hil1}) shows that $Y_{11}$ and $Y_{12}$ are analytic everywhere,
and if an entire function is $O(|z|^m)$ as $z\to\i$, then it is a polynomial
of degree at most $m$. Thus, we get from (\ref{hil3})
that $Y_{11}(z)=z^n+\cdots$ is a monic polynomial
of degree $n$, and $Y_{21}(z)$ is a polynomial of degree at most $n-1$. 
The relation (\ref{hil4}) gives that the integral of
$z^kY_{12}(z)$ over the circle $|z|=R$ is $O(R^{k-n})$ 
for all $k< n$ and hence it
tends to 0 as $R\to\i$.
By analyticity, the 
integral over the upper part of the circle can be deformed into
an integral from $R$ to $-R$ on the upper part of $\R$,
i.e., into
\[\int_{R}^{-R}x^kY_{12}^+(x)dx,\]
 and
similarly the integral over the lower part of the circle can be deformed into
an integral from $-R$ to $R$ on the lower part of $\R$, i.e., into
\[\int_{-R}^{R}x^kY_{12}^-(x)dx.\]
The first relation in (\ref{hil2}) 
implies
\[x^kY_{12}^-(x)-x^kY_{12}^+(x)=x^kY_{11}(x)w(x),\]
therefore for $k=0,1,\ldots,n-1$ we have 
\[\int_{-R}^R x^kY_{11}(x)w(x)dx=O(R^{k-n})=O(R^{-1})\to 0,\]
which implies
\[\int_{-\i}^\i x^kY_{11}(x)w(x)dx=0.\]
Thus, $Y_{11}$ is indeed the monic $n$-th orthogonal polynomial with
respect to $w$.

The other entries
can also be explicitly written in terms  
of the orthogonal polynomials $p_n$ and $p_{n-1}$:
$Y_{21}$ is a constant multiple of $p_{n-1}$,
\[Y_{12}(z)=\frac{1}{2i\pi \k_n}\int \frac{p_n(x)w(x)}{x-z}dx\]
is the Cauchy transform of $p_n(x)w(x)/\k_n$, and $Y_{22}$ is
the Cauchy transform of $Y_{21}$ ($={\rm const}\cdot p_{n-1}$).
Furthermore, $\k_n$ and the recurrence coefficients $a_n,b_n$ 
(see Section \ref{realsect}) can be expressed
in terms of the entries of $Y_1$, where $Y_1$ is the matrix defined by
\[Y(z)\left(\begin{array}{cc}
  z^{-n} & 0 \\
  0 & z^n
\end{array}\right)=:I+z^{-1}Y_1+O\left(\frac{1}{z^2}\right).\]

For details on this Riemann--Hilbert approach, see Deift \cite{Deift1}.

\subsubsection{Orthogonal polynomials with respect to inner products}
Sometimes one talks about orthogonal polynomials with respect to
an 
inner product $\langle\cdot,\cdot\rangle$ which is defined on some
linear space containing all polynomials, and orthogonality
means $\langle p_n,p_m\rangle=0$ for $ m\not=n$.
In this case the aforementioned
orthogonalization
process can be used, and with $\s_{i,j}=\langle x^i,x^j\rangle$,
the determinantal formula (\ref{det}) is still valid.

Sometimes one has an  $\langle\cdot,\cdot\rangle$ with the
standard inner product properties, except that 
positive definiteness may not hold (as an example
consider non-Hermitian
orthogonality  from  Section \ref{secnonherm}). Then
the orthogonalization process and the determinantal formulae
can still be used provided  the Gram determinants (\ref{gram})
are different from zero. If
this is not so, then 
we write
\[p_n(z)=\g_nz^n+\g_{n-1}z^{n-1}+\cdots,\]
and make sure that $p_n$ is orthogonal to all powers $z^k$, $0\le
k<n$, i.e., solve the homogeneous system of equations
\[\sum_{j=0}^n\g_j\s_{j,k}=0,\qquad k=0,\ldots,n-1,\]
for $\g_0,\g_1,\ldots,\g_n$. Since the number of unknowns is
bigger than the number of equations, there is always a non-trivial
solution, which gives rise to non-trivial orthogonal polynomials. However, then
we cannot assert any more $\g_n\not=0$, so the degree of $p_n$ may be smaller
than $n$, and there may be several choices for $p_n$. Still, in applications
where non-Hermitian orthogonality is used, these $p_n$ play the role of orthogonal
polynomials.

\subsubsection{Varying weights}
In the last 25 years orthogonal polynomials with respect to \dword{varying measures}
have played a significant role in several problems, see e.g. the sections
on exponential and Freud weights or on random matrices in Section 4.
In forming them one has a sequence of measures $\mu_n$ (generally with some
particular behavior), and for each $n$ one forms the orthogonal system
$\{p_k(\mu_n,z)\}_{k=0}^\i$. In most cases one needs the behavior of
$p_n(\mu_n,z)$ or that of $p_{n\pm k}(\mu_n,z)$ with some
fixed $k$. We mention three examples.

The first example is that of Freud weights: $W(x)=e^{-|x|^\l}$, $\l>0$.
If one substitutes $x=n^{1/\l}y$, then with $P_n(y)=p_n(W,x)$
orthogonality takes the form
\[\int P_n(y)P_m(y)e^{-n|y|^\l}dy=0,\qquad n\not=m,\]
and it turns out that this is just the right scaling, e.g. 
the zeros of $P_n$ have an asymptotic zero distribution (while
those of $p_n(W,z)$ are spreading out to infinity). Thus, studying 
orthogonal polynomials with respect to Freud weights $W$
is equivalent to studying orthogonal polynomials
with respect to the varying weights $w_n(x)=W(x)^n$, and actually,
working with $w_n$ turns out to be very natural.

For the second and third examples see multipoint Pad\'e approximation 
and random matrix theory in Section \ref{multipoint}.

\subsubsection{Matrix orthogonal polynomials}
Orthogonality of matrix polynomials (i.e., when the entries
of the fixed size matrix are polynomials of degree $n=0,1,\ldots$
and orthogonality is with respect to a matrix measure)
is a very active area which shows extreme richness compared to
the scalar case. See Section \ref{secmatrix} for a short
discussion.

\sect{The $L^2$ extremal problem}
One of the most useful tools in the study of orthogonal polynomials is
the fact that they solve the following extremal problem: minimize
the $L^2(\mu)$ norm for all monic polynomials $P_n(z)=z^n+\cdots$
of degree $n$. The minimum turns out to be  $1/\k_n^2$, 
i.e., the 
$n$-th monic orthogonal polynomial is the (unique) extremal
polynomial in
\be \inf_{P_n(z)=z^n+\cdots}\int |P_n|^2d\mu=\frac1{\k_n^2}.\label{extremal}\ee
Indeed, any $P_n$ is a linear combination $\sum_{k=1}^nc_k p_k$ with
the orthonormal polynomials $p_k$, and here
$c_n=1/\k_n$ because $P_n$ is monic, i.e., it has leading coefficient 1. 
Orthogonality gives
\[\int |P_n|^2d\mu=\sum_{k=0}^n|c_k|^2,\]
from which we can see that this is always $\ge |c_n|^2=1/\k_n^2$,
and equality occurs if and only if all the other $c_k$'s are 0.

A related extremum problem leads to the so called \dword{Christoffel functions}
associated with  $\mu$. They are defined as
\be \l_n(\mu,z)=\inf_{P_n(z)=1,\ {\rm deg}(P_n)\le n} \int |P_n|^2d\mu.
\label{christ}\ee
If we write again $P_n=\sum_{k=0}^n c_kp_k(\mu,\cdot)$, then
$P_n(z)=1$ means 
\[\sum_{k=0}^nc_kp_k(\mu,z)=1,\]
 and hence by Cauchy's inequality
\[1\le (\sum_{k=0}^n|c_k|^2) (\sum_{k=0}^n|p_k(\mu,z)|^2).\]
Therefore,
\[\int |P_n|^2d\mu=\sum_{k=0}^n|c_k|^2\ge  (\sum_{k=0}^n|p_k(\mu,z)|^2)^{-1}\]
with equality if and only if 
\[c_k=\frac{\overline{p_k(\mu,z)}} {\sum_{k=0}^n|p_k(\mu,z)|^2}.\]
Thus, we have arrived at the formula
\be\l_n(\mu,z)^{-1}=\sum_{k=0}^n|p_k(\mu,z)|^2\label{christformula}\ee
for all $z\in \C$ for the Christoffel function $\l_n(\mu,z)$.

For example, for measures $\mu$ lying on the real line
it is easy to see from this formula that
$\mu$ has a point mass at $x_0$, i.e., $\mu(\{x_0\})>0$
if and only if $\sum_kp_k(\mu,x_0)^2<\i$, and then
\[\mu(\{x_0\})=(\sum_{k=0}^\i p_k(\mu,x_0)^2)^{-1}.\]

\sect{Orthogonal polynomials on the real line}\label{realsect}
Let $\mu$ be supported on the real line. In this case orthogonalization
leads to real polynomials (i.e., all the coefficients are real). The most remarkable
property of this real case is that the
$p_n$'s obey a \dword{three-term recurrence formula} \be
xp_n(x)=a_np_{n+1}(x)+b_np_n(x)+a_{n-1}p_{n-1}(x),\label{rek}\ee
where
\be a_n=\frac{\k_n}{\k_{n+1}}>0, \qquad b_n=\int xp_n^2(x)d\mu(x)\label{rr2}\ee
are the so called \dword{recurrence coefficients}.
Indeed, if we write $xp_n(x)$ as a linear combination 
$\sum_{k=0}^{n+1} c_kp_k(z)$
with 
\[c_k:=\int xp_n(x) p_k(x)d\mu(x),\]
then all the $c_k$'s for $k<n-1$ vanish
by orthogonality:
\[c_k=\int xp_n(x) p_k(x)d\mu(x)=\int (xp_k(x))p_n(x)d\mu(x)=0\]
because $xp_k(x)$ is a polynomial of degree smaller than $n$.
Comparison of the leading coefficients on both
sides  gives that $c_{n+1}=\k_n/\k_{n+1}$, but since
$c_{n+1}$ is also the integral of $xp_n(x)p_{n+1}(x)$ against $\mu$, we 
get that 
\[c_{n-1}=\int xp_n(x)p_{n-1}(x)d\mu(x)=\frac{\k_{n-1}}{\k_n}.\]
Finally, $c_n$ is the integral given in (\ref{rr2}).

We emphasize that the three-term recurrence
is a very special property of real orthogonal polynomials, and it
is due to the fact that in this case the polynomials are real, hence
\[\int xp_n(x)\overline{p_k(x)}d\mu(x)=\int p_n(x)\overline{(xp_k(x))}d\mu(x)=0\]
for $k<n-1$. In the non-real case the two sides here are
totally different.  The three-term recurrence is missing in the
general case, and it is replaced by a different recurrence for polynomials on the
circle.

Conversely, any system of polynomials satisfying (\ref{rek}) with
real $a_n>0$, $b_n$ is an orthonormal system with respect to a
(not necessarily unique) measure on the real line (Favard's theorem).
The unicity of the measure in question is the same as the determinacy of
the moment problem, which in turn is again closely related to the behavior
of orthogonal polynomials; see Section \ref{determinacy}.

In the real case
the zeros of $p_n$ are real and simple and
 the zeros of $p_n$ and $p_{n+1}$ interlace, i.e., in between any two zeros of
$p_{n+1}$ there is a zero of $p_n$. In fact, $p_n$ must have $n$ sign changes,
for if it had only $m<n$, say at the points $y_1,\ldots,y_n\in \R$, then
it could not be orthogonal to the polynomial $q(x)=\prod_{j=1}^m (x-y_j)$
of degree $m<n$, for then $q(x)p_n(x)$ would be of constant sign. 
Let now $x_n<x_{n-1}<\ldots<x_1$ be the zeros of $p_n$,
and suppose that we already know that the zeros of $p_n$ and $p_{n-1}$ interlace,
which implies sign$(p_{n-1}(x_k))=(-1)^{k-1}$.
If we substitute $x_k$ into the recurrence formula
(\ref{rek}) then $a_n>0$ gives that $p_{n+1}(x_k)$ and 
$p_{n-1}(x_k)$ are of opposite signs at $x_k$, i.e., sign$(p_{n+1}(x_k))=(-1)^k$,
and this gives that 
the zeros of $p_n$ and $p_{n+1}$ also interlace. Thus, the interlacing property follows
by induction.

The three-term recurrence
implies for the so called \dword{reproducing kernel}
the \dword{Christoffel-Darboux formula}
\be \sum_{k=0}^np_k(x)p_k(t)=\frac{\k_n}{\k_{n+1}}
\frac{p_{n+1}(x)p_n(t)-p_n(x)p_{n+1}(t)}{x-t}.\label{cdformula}\ee
Indeed,  use the recurrence formula for $p_{n+1}$ on the right
and  $a_n=\k_n/\k_{n+1}$;
then induction gives (\ref{cdformula}). The special case
\be \l_n(\mu,x)^{-1}=\sum_{k=0}^np_k(x)^2=\frac{\k_n}{\k_{n+1}}\left(
p_{n+1}'(x)p_n(x)-p_n'(x)p_{n+1}(x)\right)\label{cdformula'}\ee
is worth mentioning.

The starting values of the recurrence (\ref{rek}) are $p_{-1}\equiv 0$,
$p_0=\left(\mu(\C))\right)^{-1/2}$. If one starts from $q_{-1}=-1$,
$q_0\equiv 0$ and uses the same recurrence (with $a_{-1}=1$)
\be xq_n(x)=a_nq_{n+1}(x)+b_nq_n(x)+a_{n-1}q_{n-1}(x),\label{rek1}\ee
then $q_n$ is of degree  $n-1$, and by Favard's theorem the different
$q_n$'s  are orthogonal
with respect to some measure. The $q_n$'s are called
\dword{orthogonal polynomials of the second kind} 
(sometimes for $p_n$ we say that they are \dword{of the first kind}).
They can also be written in the form
\[q_n(z)=\left(\mu(\C)\right)^{-1/2}\int\frac{p_n(z)-p_n(x)}{z-x}d\mu(x).\]

\sect{Classical orthogonal polynomials}\label{secclassic} These are

\begin{itemize}
  \item \dword{Jacobi polynomials} $P_n^{(\a,\b)}$, $\a,\b>-1$, orthogonal with respect to the weight
  $(1-x)^\a(1+x)^\b$ on $[-1,1]$,
  \item \dword{Laguerre polynomials} $L_n^{(\a)}$, $\a>-1$,
   with orthogonality weight $x^\a e^{-x}$
  on $[0,\i)$,
  \item \dword{Hermite polynomials} $H_n$ orthogonal with respect to $e^{-x^2}$ on 
  the real line $(-\i,\i)$.
\end{itemize}
In the literature various normalizations are used for them.

They are very special, for they possess many properties that no other orthogonal polynomial system does. In particular,
\begin{itemize}
\item they have derivatives which form again an orthogonal polynomial system, e.g.
  the derivative of $P_n^{(\a,\b)}$ is a constant multiple of $P_{n-1}^{(\a+1,\b+1)}$:
  \[(P_n^{(\a,\b)})'(x)=\frac12 (n+\a+\b+1)P_{n-1}^{(\a+1,\b+1)}(x),\]
\item they all possess a Rodrigue's type formula
  \[P_n(x)=\frac{1}{d_nw(x)}\frac{d^n}{dx^n}\{w(x)\s(x)^n\},\]
  where $w$ is the weight function and  $\s$ is a polynomial that is independent of $n$,
  for example,
  \[L_n^{(\a)}(x)=e^xx^{-\a}\frac{1}{n!}\frac{d^n}{dx^n}\left( e^{-x}x^{n+\a}\right),\]
\item they satisfy a differential-difference relation of the form
  \[\pi(x)P_n'(x)=(\a_nx+\b_n)P_n(x)+\g_nP_{n-1}(x),\]
e.g.
\[x(L_n^{(\a)})'(x)=nL_n^{(\a)}(x)-(n+\a)L_{n-1}^{(\a)}(x),\]

  \item they satisfy a non-linear equation of the form
\bean \s(x)\left(P_n(x)P_{n-1}(x)\right)'&=&
(\a_nx+\b_n)P_n(x)P_{n-1}(x)\\
&&\qquad+\g_nP_n^2(x)+\d_nP_{n-1}^2(x),\eean
with some constants $\a_n,\b_n,\g_n,\d_n$,
and $\s$ a polynomial of degree at most 2,
e.g.
\[(H_n(x)H_{n-1}(x))'=2xH_n(x)H_{n-1}(x)-H_n^2(x)+2nH_{n-1}^2(x).\]
\end{itemize}
Every one of these properties has a converse, namely if a system of orthogonal polynomials
possesses any of these properties, then it is (up to a change of variables)
one of the classical systems, see Al-Salam \cite{alsalam}. See also Bochner's result in
the next section claiming that
the classical orthogonal polynomials are
essentially the only polynomial (not just orthogonal polynomial)
systems that satisfy
a certain second order differential equation.

Classical orthogonal polynomials are also special in the sense that they possess a relatively
simple
\begin{itemize}
  \item second order differential equation, e.g.
  \[xy''+(\a+1-x)y'+ny=0\]
  for $L_n^{(\a)}$,
  \item generating function, e.g.
  \[\sum_n\frac{H_n(x)}{n!}w^n=\exp(2xw-w^2),\]
  \item integral representation, e.g.
  \[(1-x)^\a(1+x)^\b P_n^{(\a,\b)}(x)=\frac{(-1)^n}{2^{n+1}\pi i}\int(1-t)^{n+\a}(1+t)^{n+\b}(t-x)^{-n-1}dt\]
over an appropriate contour,
\end{itemize}
and these are powerful tools to study their behavior.

For all these results see Szeg\H{o} \cite{Szego}.

\sect{Where do orthogonal polynomials come from?}\label{secwhere}
In this section we mention a few selected areas where orthogonal polynomials
naturally arise.

\subsubsection{Continued fractions}
Continued fractions played an extremely important role
in the development
of several branches of
mathematics, but their significance
has unjustly diminished in modern
mathematics. A \dword{continued fraction} is of the form
$${B_1\over A_1+{B_2\over A_2+\cdots}},$$
and its \dword{$n$-th convergent} is
$${S_n\over R_n}={B_1\over A_1+{B_2\over A_2+\cdots {B_n\over A_n}}},\quad n=1,2,\ldots.$$
The value of the continued fraction is the limit of its convergents.
The denominators and numerators of the convergents satisfy the 
three-term recurrence relations
\bean
R_n&=&A_n R_{n-1}+B_nR_{n-2},\qquad R_0\equiv 1,\ R_{-1}\equiv 0\\
S_n&=&A_n S_{n-1}+B_nS_{n-2},\qquad S_0\equiv 0,\ S_{-1}\equiv 1,\eean
which immediately connects continued fractions with three-term recurrences
and hence with orthogonal polynomials.

But the connection is deeper than just this formal observation.
Many elementary functions (like $z-\sqrt{z^2-1}$) have a continued
fraction development where the $B_n$'s are constants while the
$A_n$'s are linear functions, in which case the convergents are
ratios of some orthogonal polynomials of the second and first
kind. An important example is that of Cauchy transforms of
measures $\mu$ with compact support on the real line (so called
\dword{Markov functions}): \be f(z)=\int{d\mu(x)\over x-z}=-{\a_0\over
z}-{\a_1\over z^2}-\ldots.\label{markov}\ee
 The coefficients $\a_j$ in the development of (\ref{markov}) are the moments
$$\a_j=\int x^jd\mu(x),\quad j=0,1,\ldots$$
of the measure $\mu.$ The continued fraction development
$$f(z)\sim {B_1\over z-A_1+{B_2\over z-A_2+\cdots}}$$
of $f$ at infinity converges locally uniformly outside the smallest interval
that contains the support of $\mu$ (A. Markov's theorem).

As has been mentioned, the
numerators $S_n(z)$ and the denominators $R_n(z)$ of the $n$-th convergents
$${S_n(z)\over R_n(z)}={B_1\over z-A_1+{B_2\over z-A_2+\cdots {B_n\over z-A_n}}},\quad n=1,2,\ldots$$
satisfy the recurrence relations
\bea
R_n(z)&=&(z-A_n)R_{n-1}(z)+B_nR_{n-2}(z),\quad R_0\equiv 1,\ R_{-1}\equiv 0\label{rekk}\\
S_n(z)&=&(z-A_n)S_{n-1}(z)+B_nS_{n-2}(z),\quad S_0\equiv 0,\ S_{-1}\equiv 1.\nonumber\eea
These are precisely the recurrence formulae for the monic orthogonal polynomials of the
first and second kind with respect to $\mu$, hence the $n$-th convergent
is $cq_n(z)/p_n(z)$ with $c=\mu(\C)^{1/2}$.

See Szeg\H{o} \cite[pp. 54--57]{Szego} as well as Kruschev 
\cite{Kruschev} and the
numerous references there.

\subsubsection{Pad\'e approximation and rational interpolation}\label{multipoint}
One easily gets from the recurrence relations (\ref{rekk})
that
\[\frac{S_m(z)}{R_m(z)}-\frac{S_{m+1}(z)}{R_{m+1}(z)}=(-1)^n\frac{B_1B_2\cdots B_{n+1}}
{R_n(z)R_{n+1}(z)},\]
and summation of these for $m=n,n+1,\ldots$ yields that
\[\frac{S_n(z)}{R_n(z)}=\sum_{k=0}^{2n} \frac{-\a_k}{z^{k+1}}+O(z^{-2n-1}),\]
i.e., with the preceding notation the rational function 
\[S_n(z)/R_n(z)=cq_n(z)/p_n(z)\quad\mbox{with $c=\mu(\C)^{1/2}$}\]
 of numerator degree $n-1$ and
of denominator degree $n$ interpolates $f(z)$ at infinity to order $2n$.
This is the analogue (called $[n-1/n]$ \dword{Pad\'e approximant}) of
the $n$-th Taylor polynomial
(which interpolates the function to order $n$) for rational functions.
The advantage of Pad\'e approximation over Taylor polynomials lies in the fact
that the poles of Pad\'e approximants may imitate the singularities
of the function in question, while Taylor polynomials are good only up
to the first singularity. The error in $[n-1/n]$ Pad\'e approximation
has the form
\[f(z)-c{q_n(z)\over p_n(z)}={1\over p_n^2(z)}\int{p_n^2(x)\over x-z}d\mu(x).\]

Orthogonal polynomials appear in more general rational interpolation
(called \dword{multipoint Pad\'e approximation}) to Markov functions,
see e.g. Stahl--Totik \cite[Sec. 6.1]{stahltotik}. For every $n$  select a set $A_n=\{x_{0,n},\ldots,x_{2n,n}\}$ of $2n+1$ interpolation points from 
$\overline\C\setminus I$ where $I$ is the smallest
interval that contains the support of $\mu$. The points need not
be distinct, but we assume that $A_n$ is symmetric with respect to the
real line. Put
\[ \omega_n(z):=\prod^{2n}_{\scriptstyle j=0\atop\scriptstyle  x_{jn}\not=\infty}(z-x_{jn}).\]
The degree of $\omega_n$ is equal the number of finite points in $A_n.$ 
By  $r_n(z)=u_n(z)/Q_n(z)$ we denote a rational function of numerator and denominator
degree at most $n$ that interpolates the function $f$ at the $2n+1$ points of the set $A_n=\{x_{0,n},\ldots,x_{2n,n}\}$ in the sense that
\[{f(z)-r_n(f,A_n;z)\over \omega_n(z)}=O(z^{-(2n+1)})\quad {\rm as}\quad |z|\to\infty;\]
the expression on the left 
is bounded at every finite point of $A_n$, and at infinity it has the indicated behavior.
Now for Markov functions  this rational interpolant uniquely  exists,  
$Q_n$ is the $n$-th orthogonal polynomial with respect to the varying weight
$d\mu(x)/\o_n(x)$, and the remainder term of the interpolation has the representation
\[(f-r_n(f,A_n;\cdot))(z)={\omega_n(z)\over
 Q^2_n(z)}\int{Q_n^2(x)\over\omega_n(x) (x-z)}d\mu(x)\]
for all $z$ outside the support of $\mu$.
Thus, the rate of convergence of the rational interpolants
is intimately connected with the behavior of the orthogonal polynomials
with respect to  the varying weight $d\mu(x)/\o_n(x)$.

\subsubsection{Moment problem}\label{determinacy}
The \dword{moments} of a measure $\mu$, $\mu(\C)=1$,
 supported on the real line, are
\[\a_n=\int x^nd\mu(x), \qquad n=0,1,\ldots.\]
The Hamburger moment problem is to determine if a sequence $\{\a_n\}$
(with normalization $\a_0=1$) of real
numbers is the moment sequence of a measure with infinite support,
and if this measure is unique (the
Stieltjes moment problem asks the same, but for measures on $[0,\i)$).
The existence is easy: $\{\a_n\}$ are the moments of some measure supported
on $\R$ if and only if all the Hankel determinants $|\a_{i+j}|_{i,j=0}^m$, $m=0,1,\ldots$,
are positive. The unicity (usually called determinacy) depends on the behavior of the
orthogonal polynomials (\ref{det}) defined from the moments $\s_{i,j}=\a_{i+j}$. In
fact, there are different measures with the same moments $\a_j$ if and only if
there is a non-real $z_0$ with $\sum_n|p_n(z_0)|^2<\i$, which in turn is
equivalent to
$\sum_n|p_n(z)|^2<\i$ for all $z\in \C$. Furthermore, the Cauchy transforms of
all solutions $\nu$ of the moment problem have the
parametric form
\[\int\frac{d\nu(x)}{z-x}=\frac{C(z) F(z)+A(z)}{D(z)F(z)+B(z)},\]
where $F$ is an arbitrary analytic function (the parameter)
mapping the upper half plane $\C_+$ into $\overline
{\C}_+\cup\{\i\}$, and $A,B, C$ and $D$ have explicit
representations in terms of the first and second kind orthogonal
polynomials $p_n$ and $q_n$: \bean \begin{array}{ll}A(z)=z\sum_n
q_n(0)q_n(z); & B(z)=-1+z\sum_n q_n(0)p_n(z);
\\[10pt]
C(z)=1+z\sum_n p_n(0)q_n(z); & D(z)=z\sum_n p_n(0)p_n(z).
\end{array}\eean

For all these results see Akhiezer \cite{akhiezer},
and for an operator theoretic approach to the moment
problem see Simon \cite{Simonmoment} (in particular, Theorems
3 and 4.14).

\subsubsection{Jacobi matrices and spectral theory of self-adjoint operators}
Tridiagonal, so called \dword{Jacobi matrices}
\[J= \left(\begin{array}{ccccc}
  b_0 & a_0 & 0 & 0 &\cdots  \\
  a_0 & b_1 & a_1 & 0 & \cdots \\
  0 & a_1 & b_2 & a_2 & \cdots \\
  0 & 0 & a_2 & b_2 &\cdots \\
  \vdots&\vdots&\vdots&\vdots&\ddots
\end{array}\right)\]
with bounded $a_n>0$ and bounded real $b_n$
 define a bounded self-adjoint operator $J$
in  $l_2$, a so called \dword{Jacobi operator}. These are the discrete analogues of
second order linear differential operators of  Schr\"odinger type on the half line.
Every bounded
self adjoint operator with a cyclic vector is a Jacobi operator in an appropriate
base.

The formal eigen-equation $J\pi=\l\pi$  is equivalent to the three-term recurrence
\[a_{n-1}\pi_{n-1}+b_n\pi_n+a_n\pi_{n+1}=\l \pi_n, \qquad n=1,2,\ldots,\]
\[b_0\pi_0+a_0\pi_1=\l\pi_0, \qquad \pi_0=1.\]
Thus,  $\pi_n(\l)$ is of degree $n$ in $\l$.

By the spectral theorem,
$J$, as a self-adjoint operator having a cyclic vector ($(1,0,0,\ldots)$),
is unitarily equivalent to multiplication by $x$ in some $L^2(\mu)$
with some probability measure
$\mu$ having compact support on the real line. This $\mu$ is called
the spectral measure associated with $J$ (and with its spectrum).
 More precisely,
if $p_n(x)=p_n(\mu,x)$ are the orthonormal polynomials with respect to $\mu$, and
$U$ maps the unit vector $e_n=(0,\ldots,0,1,0,\ldots)$ into $p_n$,
then $U$ can be extended to a unitary operator from $l_2$ onto
$L^2(\mu)$, and if $Sf(x)=xf(x)$ is the multiplication operator
by $x$ in $L^2(\mu)$, then $J=U^{-1}SU$. The recurrence coefficients
for $p_n(\mu,x)$ are precisely the $a_n$'s and $b_n$'s from the
Jacobi matrix, i.e., $p_n(x)=c\pi_n(x)$ with some fixed constant $c$.
These show that Jacobi operators are
equivalent to multiplication by $x$ in
$L^2(\mu)$ spaces if the particular basis
$\{p_n(\mu,\cdot)\}$ is used (see e.g.
Deift \cite[Ch. 2]{Deift1}).

The truncated $n\times n$ matrix
\[J_n= \left(\begin{array}{ccccc}
  b_0 & a_0 & 0 & 0 &\cdots  \\
  a_0 & b_1 & a_1 & 0 & \cdots \\
  0 & a_1 & b_2 & a_2 & \cdots \\
  \vdots&\vdots&\vdots&\ddots&\cdots\\
  0&0&0&a_{n-2}&b_{n-1}
\end{array}\right)\]
has $n$ real and distinct eigenvalues, which turn out to be the zeros
of $p_n$, i.e., the monic polynomial $p_n(z)/\k_n$ is the characteristic
polynomial of $J_n$.

\subsubsection{Quadrature}
For a measure $\mu$,
an \dword{$n$-point quadrature} (rule)
is a sequence of $n$ points $x_{n,1},\ldots,x_{n,n}$ and
a sequence of associated numbers $\l_{n,1},\ldots,\l_{n,n}$. It is expected that
\[\int fd\mu\sim \sum_{k=1}^n \l_{n,k} f(x_{n,k})\]
in some sense for as large a class of functions as possible.
Often the accuracy of the quadrature is measured by its
exactness, which is defined as the largest $m$ such that
the quadrature is exact for all polynomials of degree at
most $m$, i.e., $m$ is the largest number with the property
that
\[\int x^jd\mu(x)=\sum_{k=1}^n \l_{n,k} x_{n,k}^j \qquad \mbox{for all $0\le j\le m$}.\]
For $\mu$ with support on the real line and for quadrature based on
$n$ points this exactness $m$ cannot be larger than
$2n-1$, and this optimal value $2n-1$ is attained
if and only if $x_{n,1},\ldots,x_{n,n}$ are the
zeros of the orthonormal polynomial $p_n(\mu,x)$ corresponding to
$\mu$, and the so called Cotes numbers $\l_{n,k}$ are chosen to be
\[\l_{n,k}=\l_n(\mu,x_{n,k})=\left(\sum_{j=0}^np_j(\mu,x_{n,k})^2\right)^{-1},\]
where $\l_n(\mu,z)$ is the Christoffel function (\ref{christ}) associated
with $\mu$.

See Szeg\H{o} \cite[Ch. XV]{Szego}.

\subsubsection{Random matrices}
Some statistical-mechanical models in quantum systems
use random matrices.
Let ${\cal H}_n$ be the set of all $n\times n$ Hermitian
matrices $M=(m_{i,j})_{i,j=1}^n$, and let there be given a
probability distribution on ${\cal H}_n$ of the form
\[P_n(M)dM=D_n^{-1}\exp(-n{\rm Tr}\{V(M)\})dM,\]
where $V(\l)$, $\l\in \R$, is a real-valued function
that increases sufficiently fast at infinity
 (typically an even polynomial in quantum field theory applications),
Tr$\{H\}$ denotes the trace of the matrix $H$,
\[dM=\prod_{k=1}^ndm_{k,k}\prod_{k<j}d\, \Re m_{k,j}\;d\, \Im
m_{k,j}\] is the ``Lebesgue'' measure for Hermitian matrices,
and $D_n$ is a normalizing constant so that the total integral of
$P_n(M)dM$ is one.

Every matrix $M\in {\cal H}_n$ has $n$ real eigenvalues
which carry physical information on the system when it
is in the state described by $M$.  The quantity
\[N_n(\D)=\frac{\#\{\mbox{eigenvalues\ in\ $\D$}\}}{n}\]
is  the random variable that
equals the normalized number of eigenvalues in the interval
$\D$.
This model is known as the  unitary
ensemble associated with $V$.

Let $p_j(w^n,x)$ be the  orthonormal
polynomials
 with respect to
 the varying weight $w^{n}(x)$, $w(x)=\exp(-
V(x))$. Then the joint probability density of the
eigenvalues can be written in the form
\[d_n \left|p_{i-
1}(w^n,\l_j)w^{n/2}(\l_j)\right|^2_{1\le i,j\le n},\] where $d_n$
is a normalizing constant built up from the leading coefficients
of the $p_j(w^n,\cdot)$. With the so called weighted reproducing
kernel
\[K_n(t,s)=\sum_{j=0}^{n-
1}p_j(w^n,t)w^{n/2}(t)\,p_j(w^n,s)w^{n/2}(s),\]
it can also be written in the form
\[\frac{1}{n!}\left|K_n(\l_i,\l_j)\right|_{1\le i,j\le n}.\]
In particular, for the expected number of eigenvalues
in an interval $\D$ we have
\[EN_n(\D)=\int_\D\frac{K_n(\l,\l)}{n}\;d\l,\]
where $1/K_n(\l,\l)$ is
known in the theory of orthogonal
polynomials as the $n$-th (weighted) Christoffel
function  associated with the
weight $w^n$, while the limit of the left hand side
(as $n\to\i$)
is known as the density of states.

See, e.g., Mehta \cite{mehta} and Pastur--Figotin \cite{Pastur}.

\section{Some questions leading to classical orthogonal polynomials}
There are almost an infinite number of problems where classical orthogonal
polynomials emerge. Let us just mention a few.
\bigskip

\subsubsection{Electrostatics}
\smallskip

Put at $1$ and $-1$ two positive charges
$p$ and $q$, and with these fixed charges put $n$ positive unit charges
on $[-1,1]$ at the points $x_1,\ldots,x_n$. On the plane the Coulomb force
is proportional with the reciprocal of the distance, and so a charge
generates a logarithmic potential field. Therefore,
the mutual energy of all these
charges is
\[I(x_1,\ldots,x_n)=p\sum_{j=1}^n\log\frac{1}{|1-x_j|}+
    q\sum_{j=1}^n\log\frac{1}{|1+x_j|}+
    \sum_{i<j}\log\frac{1}{|x_i-x_j|},\]
and the equilibrium problem asks for finding $x_1,\ldots,x_n$ for which
this energy is minimal. The unique minimum occurs (see Szeg\H{o} 
\cite[Section 6.7]{Szego})
for the
zeros of the Jacobi polynomial $P_n^{(2p-1,2q-1)}$ orthogonal with respect to
the weight $(1-x)^{2p-1}(1+x)^{2q-1}$.

There is a similar characterization of the zeros of Laguerre and Hermite polynomials,
and even of more general orthogonal polynomials (for the latter see
Ismail \cite[Section 3.5]{Ismail}).
\bigskip

\subsubsection{Polynomial solutions of eigenvalue problems}
\smallskip

Consider the eigenvalue problem
\[f(x)\frac{d^2}{dx^2}y(x)+g(x)\frac{d}{dx}y(x)+h(x)y(x)=\l y(x),\]
where $f,g,h$ are fixed polynomials and $\l$ is a free constant, and it
is required that this have a polynomial solution of exact
degree $n$ for all $n=0,1,\ldots$, for which the corresponding $\l$ and $y(x)$
will be denoted by $\l_n$ and $y_n(x)$, respectively.
Bochner's theorem  from \cite{bochner}
states that, except for some trivial solutions of the form $y(x)=ax^n+bx^m$
and for some
polynomials related to Bessel functions, the only solutions are (in all of them
we can take $h(x)=0$)
\begin{itemize}
\item Jacobi polynomials $P_n^{(\a,\b)}$ ($f(x)=1-x^2$, $g(x)=\b-\a-x(\a+\b+2)$,
$\l_n=-n(n+\a+\b+1)$)

\item Laguerre polynomials $L_n^{(\a)}$ ($f(x)=x$, $g(x)=1+\a-x$, $\l_n=-n$) and

\item Hermite polynomials $H_n(x)$ ($f(x)=1$, $g(x)=-2x$, $\l_n=-2n$).
\end{itemize}
\bigskip

\subsubsection{Harmonic analysis on spheres and balls}
\smallskip

Harmonic analysis on spheres and balls in $\R^d$ is based on spherical
harmonics, i.e., harmonic homogeneous polynomials. In this theory, special
Jacobi polynomials, the so called \dword{ultraspherical} or \dword{Gegenbauer polynomials}
$P_n^{(\a)}$, play a fundamental role -- they are orthogonal with respect to
the weight $(1-x^2)^{\a-1/2}$.

Let $S^{d-1}$ be the unit sphere in $\R^d$ and let ${\cal H}_n^d$ be the restriction
to $S^{d-1}$ of all harmonic polynomials $Q(x_1,\ldots,x_n)$
of $d$ variables that are
homogeneous of degree $n$, i.e.,
\[\sum_{k=1}^n\frac{\partial^2}{\partial x_k^2}Q=0,\qquad
Q(\l x_1,\ldots,\l x_n)=\l^nQ(x_1,\ldots,x_n),\quad \l>0.\]
The dimension of ${\cal H}_n^d$ is
\[{n+d-1 \choose d-1}-{n+d-3\choose d-1},\]
and an orthogonal basis in it
can be produced as follows. With $\r=x_{d-1}^2+x_d^2$
let $g_{s,0}=\r^s P^{(0)}_s(x_{d-1}/\r)$ and $g_{s,1}=x_d \r^s P^{(1)}_s(x_{d-1}/\r)$.
With $n_d=0$ or $n_d=1$ consider all multiindices
${\bf n}=(n_1,n_2,\ldots,n_d)$ such that $n_1+\cdots+n_d=n$, and
if for such a multiindex we define the function 
$Y_{\bf n}(x_1,\ldots,x_d)$
as
\[g_{n_{d-1},n_d}\prod_{j=1}^{d-2}\left(
(x_j^2+\cdots+x_d^2)^{n_j}P_{n_j}^{(\l_j)}(x_j(x_j^2+\cdots+x_d^2)^{-1/2})\right),\]
then these $Y_{\bf n}$ constitute an orthogonal basis in ${\cal H}_n^d$
(see e.g.  Dunkl--Xu \cite[p. 35]{xu}).

If $\underline x=(x_1,\ldots,x_n)$ and 
$\langle \underline x,\underline y\rangle=\sum_kx_ky_k$
is the inner product in $\R^d$, then the reproducing kernel for
these spherical polynomials is $P_n^{(d-2)/2}(\langle \underline x,\underline y\rangle)$
in the sense that for all $Q\in {\cal H}_n^d$ and for all $x\in S^{d-1}$ we have
\[c_{n,d}\int P_n^{(d-2)/2}(\langle \underline x,\underline y\rangle) Q(\underline y)d\s(\underline y)=Q(\underline x),\]
where integration is with respect to surface area, and $c_{n,d}$ is an
explicit normalizing constant (see e.g. Dunkl--Xu \cite[p. 37]{xu}). 

As a result, Gegenbauer polynomials are all over the theory
of spherical harmonics,
as well as in the corresponding theory for the unit ball.

\bigskip

\subsubsection{Approximation theory}
\smallskip

In the literature, expansions of functions into classical orthogonal
polynomial series are  second only to trigonometric expansions, and
numerous works have been devoted to their convergence and approximation properties,
see e.g. Szeg\H{o} \cite[Ch. XIII]{Szego}.

The \dword{Chebyshev polynomials}
$\cos(n \arccos x)$
are orthogonal on $[-1,1]$ with respect to
the weight $w(x)=(1-x^2)^{-1/2}$.
These directly
correspond to trigonometric functions, and expansions into them have
virtually the same properties as trigonometric Fourier expansions. But there
are many other aspects of approximation where Chebyshev polynomials appear.
If one considers, for example, the best approximation on
$[-1,1]$ of $x^n$ in the uniform norm
by polynomials $P_{n-1}(x)$ of smaller degree then the smallest
error appears when $x^n-P_{n-1}(x)=2^{1-n}\cos(n\arccos x)$ is the
monic $n$-th Chebyshev polynomial. Actually, monic
Chebyshev polynomials
minimize all $L^p(w)$, $p>0$, norms among monic polynomials of a given
degree. 

As we have seen in  (\ref{extremal}), the monic orthogonal 
polynomials $p_n(\mu)/\k_n$ are the solutions to the
extremal problem
\be \int |P_n|^2d\mu\to {\rm min},\label{min}\ee
where the minimum is taken for all monic polynomials of degree $n$.
This extremal property makes orthogonal polynomials, in particular
Chebyshev polynomials, indispensable tools in approximation theory.

Lagrange interpolation and its various generalizations like
Hermite\--Fej\'er or Hermite interpolation etc. is mostly done on the zeros of
some orthogonal polynomials. In fact, these nodes are often close to optimal
in the sense that the Lebesgue constant increases at the optimal rate.
In many cases interpolation on zeros of orthogonal polynomials
has special properties due to
explicitly calculable expressions. Recall e.g. Fej\'er's result
that if $P_{2n-1}$ is the unique polynomial of degree
at most $2n-1$ that interpolates
a continuous function $f$ at the nodes of the
$n$-th Chebyshev polynomial
and that has zero derivative at each of these nodes,
then $P_{2n-1}$ uniformly converges to $f$ on $[-1,1]$ as $n\to\i$.
For the role of orthogonal polynomials in interpolation see
the books Szabados--V\'ertesi \cite{szabados} and Mastroianni--Milovanovic
\cite{mastroianni}.

\sect{Heuristics}\label{secheuristics}

In this section we do not state precise results. We just want to indicate
some heuristics on the behavior of orthogonal polynomials. For the concepts below,
as well as for a more precise form of some of the heuristics see the
following sections, in particular
Section \ref{secgeneralop}.

As we have seen, the monic orthogonal polynomials $p_n(\mu)/\k_n$
minimize the $L^2(\mu)$ norm in (\ref{min}). Therefore, these polynomials
try to be small where the measure is large, e.g. one expects the zeros
to cluster at the support $S(\mu)$ of $\mu$. The example of arc measure on the unit
circle, for which the orthogonal polynomials are $z^n$, shows
however, that this is not
true (due to the fact that the complement of the support
is not connected). The statement is true when the
support lies on $\R$ or on some systems of arcs,
and also in the general case when instead of the support one considers
the polynomial convex hull
of the support of $\mu$ (for the definition, see the next
section): on any compact set outside the polynomial
convex hull there can only be a fixed number of zeros of $p_n(\mu)$ for every
$n$. When the complement of $S(\mu)$ is connected and $S(\mu)$ has no
interior, then
the distribution of the zeros shows a remarkable universality and
indifference
to the size of $\mu$. In many situations the distribution of
the zeros is the equilibrium distribution of the support $S(\mu)$.
When $S(\mu)=[-1,1]$, this means that under very weak assumptions
the zero distribution is always the arcsine distribution $dx/\pi\sqrt{1-x^2}$.

The $L^2(\mu)$ minimality of $p_n(\mu)/\k_n$ in the sense of (\ref{min})
is something like
minimality in the $L^\i$ norm on $S(\mu)$.
Therefore, $p_n(\mu)/\k_n$ should behave like
the monic polynomial $T_n$ minimizing
the $L^\i$ norm on $S(\mu)$ (so called Chebyshev polynomials for
$S(\mu)$). Since
\[\frac{1}{n}\log |T_n(z)|=\int \log|z-t|d\nu_n(t)\]
where $\nu_n$ has mass $1/n$ at each zero of $T_n$, in the limit the
behavior should be like
\be U^\nu (z)=\int\log|z-t|d\nu(t),\label{unu}\ee
where $\nu$ is the probability measure
on $S(\mu)$ for which the maximum of $U^\nu$ on $S(\mu)$ is as small as possible
 (this is the so called equilibrium measure of $S(\mu)$). More generally,
if $d\nu=d\nu_n=w^n(x)dx$ is a varying weight in the specified way,
then the same reasoning leads to a behavior like (\ref{unu}), but now
$\nu$ is a measure for which the supremum of $U^\nu(z)+\log w(z)$ is as small
as possible (weighted equilibrium measure).

Universal behavior can also be seen for the polynomials
themselves. Usually they obey
\be\frac1n\log|p_n(\mu,z)|\to g_{\C\setminus S(\mu)}(z,\i),
\qquad z\not\in S(\mu),\label{greenbeh}\ee
where $g_{\C\setm S(\mu)}(z,\i)$ is the Green
 function with pole at infinity
associated with the complement of the support.
When the unbounded component of the
complement of $S(\mu)$ is simply connected,
then  in that component often there is a finer
asymptotic behavior of $p_n(\mu)$ of the form
\be p_n(z)\sim d_ng_\mu(z)\F(z)^n,\qquad z\not\in S(\mu),\label{strong}\ee
where $\F$ is the mapping function that maps
$\C\setminus S(\mu)$ conformally onto
the outside of the unit disk, and $g_\mu$ is a function
(might be called generalized Szeg\H{o} function) that depends
on $\mu$. Such a fine asymptotic is restricted to
the simply connected case, see e.g. Section \ref{secreal}.

Asymptotics of orthogonal polynomials have a hierarchy, and
the different types of asymptotics usually require the measure to
be sufficiently strong with different degree on its support.
Consider first the case of
compact support $S(\mu)$. The weakest
is \dword{$n$-th root asymptotics} stating the behavior
(\ref{greenbeh}) for  $|p_n(\mu,z)|^{1/n}$
outside the support of the measure. It is mostly equivalent to
a corresponding distribution of the zeros, as well as asymptotical
minimal behavior of $\k_n^{1/n}$. It holds under very weak
assumptions on the measure, roughly stating that the logarithmic capacity
of the points where $\mu'>0$ (derivative with respect to equilibrium measure),
be the same as the capacity of $S(\mu)$.
\dword{Ratio asymptotics}, i.e., asymptotic
behavior of $p_{n+1}(\mu,z)/p_n(\mu,z)$, is stronger, and is equivalent
with asymptotics for the ratio $\k_{n+1}/\k_n$ of consecutive
leading coefficients. It can only hold when $\C\setm S(\mu)$
(more precisely its unbounded component)
is simply connected, and in this case it is enough that $\mu'>0$ almost
everywhere with respect to the equilibrium measure of the support of
$\mu$ (see Section \ref{secreal}).
Finally, \dword{strong asymptotics} of the form
(\ref{strong}) needs roughly that $\log \mu'$ be integrable
(Szeg\H{o} condition, see Section \ref{secreal}).

All these are outside the support. On the support the orthogonal
polynomials are of oscillatory behavior, and in the real case under
smoothness assumptions on the measure often a so called Plancherel-Rotach type
asymptotic formula
\[p_n(\mu,x)\sim d_ng(x)\sin(nh(x)+H(x))\]
holds, where $g,h,H$ are  fixed functions. Here $h(x)$  is directly
linked with the zeros, $h'/\pi$ is precisely the distribution of
the zeros.
When $S(\mu)=[-1,1]$ and the measure is smooth, then
$h(x)=\arccos x$.

When $S(\mu)$ is not of compact support (like Laguerre, Hermite
or Freud weights), then usually the zeros are spreading out, and one
has to rescale them to $[-1,1]$ (or to $[0,1]$) to get a distribution, which is
mostly NOT the arcsine distribution. In a similar fashion, various
asymptotics hold  for the polynomials only after a corresponding rescaling.

The Christoffel function (\ref{christ}) or, what is the same, the square sum
\[\frac{1}{\l_n(\mu,z)}=\sum_{k=0}^n |p_k(\mu,z)|^2,\]
behaves much more regularly than the orthogonal polynomials.
Outside the support the behavior of $\l_n$ is what one
gets from the heuristics above on the polynomials (after square summation).
On $S(\mu)$ the typical behavior of $\l_n(\mu,z)$
is like $\mu(\D_n(z))$, where $\D_n(z)$ is the disk about
$z$ with equilibrium measure $1/n$ (equilibrium measure of the
support $S(\mu)$). In particular, $n\l_n(\mu,z)$ 
tends pointwise to the Radon-Nikodym derivative of $\mu$ with respect to
the equilibrium measure. As a rule of thumb the
estimate $|p_n(\mu,x)|^2\le C/n\l_n(\mu,x)$ holds in many cases.

If $f\in L^2(\mu)$, its Fourier expansion into $\{p_n(\mu,\cdot)\}$
is
\[f(x)\sim \sum_{k=0}^\i c_kp_k(x),\qquad c_k=\int f\overline{p_k}d\mu.\]
The $n$-th partial sum has the closed form
\[\int f(t)K_n(x,t)d\mu(t),\qquad K_n(x,t)=\sum_{k=0}^n p_k(x)p_k(t).\]
In the real case  for the reproducing kernel $K_n(x,t)$
we have the Christoffel-Darboux formula
\[K_n(x,t)=\frac{\k_n}{\k_{n+1}}
\frac{p_{n+1}(x)p_n(t)-p_n(x)p_{n+1}(t)}{x-t},\]
which suggests a singular integral-type behavior for the partial
sums. In general, Fourier expansions into orthogonal polynomials
are sensitive to the weight (recall e.g. Pollard's theorem
that Legendre expansions are bounded in $L^p[-1,1]$ only for
$4/3<p<4$), but sometimes convergence properties are
equivalent to those of a related trigonometric Fourier series
(so called transplantation theorems, 
see e.g. Askey \cite{Askey1}).

\sect{General orthogonal polynomials}\label{secgeneralop}
In this section $\mu$ always has compact support $S(\mu)$. For
all the results below see Stahl--Totik \cite{stahltotik} and the references
therein.

\subsubsection{Lower and upper bounds}
The \dword{energy} $V(K)$
of a compact
set $K$ is defined as the infimum of
\be I(\nu)= \int\int\log\frac1{|x-t|}d\nu(x)d\nu(t)\label{energy}\ee
where the infimum is taken for all positive Borel measures on
$K$ with total mass 1. The \dword{logarithmic capacity} is
then $\c(K)=e^{-V(K)}$. For the leading coefficients $\k_n$ of
the orthonormal polynomials $p_n(\mu)$ we have
\be \frac1{\c(S(\mu))}\le \liminf_{n\to\i} \k_n^{1/n}.\label{lower}\ee
To get an upper bound we need the concept of carrier: a Borel set 
$E$ is a 
\dword{carrier} for $\mu$ if $\mu(\C\setm E)=0$. 
The capacity of a Borel set is the 
supremum of the capacities of its compact subsets, and the 
\dword{minimal carrier capacity} $c_\mu$ associated with $\mu$ 
is the infimum of the capacities
of all carriers. With this
\be \limsup_{n\to\i}\k_n^{1/n}\le \frac{1}{c_\mu}.\label{cmu}\ee

When $\c(K)$ is positive, then there is a unique measure $\nu=\o_K$
minimizing the energy in (\ref{energy}), and this measure is called the
\dword{equilibrium measure} of $K$. \dword{Green's function} 
$g_{\C\setm K}(z,\i)$ with pole at
infinity
of $\C\setm K$ can then be defined as
\be g_{\C\setm K}(z,\i)=\log\frac{1}{\c(K)}-\int\log\frac{1}{|z-t|}d\o_K(t).
\label{greenrep}\ee
We have for all $\mu$ (with $\c(S(\mu))>0$) the estimate
\be \liminf_{n\to\i}\frac1{n}\log |p_n(\mu,z)|^{1/n}\ge g_{\C\setm S(\mu)}(z,\i)
\label{lest}\ee
locally uniformly outside the convex hull of $S(\mu)$, while in the convex
hull but outside the so called polynomial convex hull Pc$(S(\mu))$ 
(for the definition see below) (\ref{lest}) is true \dword{quasi-everywhere} (i.e.,
with the exception of a set of zero capacity). The same is true on
the \dword{outer boundary} of $S(\mu)$, which is defined as the boundary
$\partial \O$ of the unbounded component $\O$
of the complement $\C\setm S(\mu)$, namely for quasi-every $z\in \partial \O$
\[\liminf_{n\to\i}|p_n(\mu,z)|^{1/n}\ge 1.\]
The \dword{minimal carrier Green function} $g_\mu(z,\i)$ 
is the supremum for all carriers $E$ of the Green function
of $\C\setm E$, where the latter is defined as the infimum of 
$g_{\C\setm K}$ for all compact subsets $K$ of $E$. With this,
\be \limsup_{n\to\i}\frac1{n}\log |p_n(\mu,z)|^{1/n}\le g_\mu(z,\i)
\label{lest1}\ee
locally uniformly on the whole plane. 

 All these estimates are
sharp.

When the bounds in (\ref{lower}) and (\ref{cmu}) coincide we have convergence
for $\k_n^{1/n}$, and these bounds coincide precisely when
the bounds in (\ref{lest}) and (\ref{lest1}) do so.

\subsubsection{Zeros}
The zeros of $p_n(\mu)$ always lie in the convex hull of the support
$S(\mu)$ of the measure $\mu$. 
This is a consequence of the $L^2$ extremal property (\ref{extremal})
of orthogonal polynomials. In fact, if there was a zero $z_0$ of $p_n(\mu,z)$
outside the convex hull of the support, then we could move that zero
towards the convex hull (along a line segment that is perpendicular
to a line separating $z_0$ and $S(\mu)$).
During this move the absolute value
of the polynomial decreases at all points of $S(\mu)$ and hence
so does its $L^2(\mu)$ norm, but that is impossible
by (\ref{extremal}). 

To say somewhat more on the location of zeros we need the concept
of the polynomial convex hull. When
$\O$ is the unbounded component of the complement $\C\setminus
S(\mu)$, then Pc$(S(\mu))=\C\setm \O$ is called the 
\dword{polynomial convex hull} of $S(\mu)$ (it is the union of $S(\mu)$ with all the
``holes'' in it, i.e., with the bounded components of $\C\setm
S(\mu)$). Now the zeros cluster on Pc$(S(\mu))$ in the sense that
for any compact subset $K$ of $\O$ there is a number $N_K$
independent of $n$, such that $p_n(\mu)$ can have at most $N_K$
zeros in $K$. The proof of this is based on the following
lemma: {\it Let $V,S\subseteq\C$ be two compact sets.
If $V$ and ${\rm Pc}(S)$ are disjoint, then there exist $a<1$ and $m\in\N$ such that for arbitrary $m$ points $x_1,\ldots,x_m\in V$ there exist $m$ points 
$y_1,\ldots,y_m\in\C$ for which the rational function
\be r_m(z):=\prod^m_{j=1}{z-y_j\over z-x_j}\label{1.3.12}\ee
has on $S$ a sup-norm satisfying
\be \|r_m\|_S\le a.\label{1.3.13}\ee}
 Taking this for granted,
 assume that $V$ is a compact set contained in $\Omega.$ We apply the
lemma with $S=S(\mu)$, and let $a<1$ and $m\in\N$ be the numbers 
in the lemma. 
Let us assume that $p_n(\mu;z)$ has at least $m$ zeros $x_1,\ldots,x_m$ on $V.$ By 
the lemma there exist $m$ points 
$y_1,\ldots,y_m\in\C$ such that the rational function $r_m$ defined as in (\ref{1.3.12}) by the points $x_1,\ldots,x_m$ 
and $y_1,\ldots,y_m$ satisfies the inequality (\ref{1.3.13}). 
With $r_m$ we define the modified monic polynomial
$$q_n(z):=r_m(z)p_n(\mu;z)=z^n+\ldots,$$
For the $L^2(\mu)$ norm of this polynomial we have the estimate
$$\|q_n\|_{L^2(\mu)}\le\|r_m\|_{S(\mu)}\|p_n(\mu;\cdot)\|_{L^2(\mu)}<\|p_n(\mu;\cdot)\|_{L^2(\mu)},$$
which contradicts the minimality (\ref{extremal}) of the monic orthogonal polynomial $p_n(\mu;z).$ Hence, we have proved that $p_n(\mu;z)$ has at most $m-1$ zeros on $V$,
as was stated.

What we have said about the zeros can be sharpened
for measures on the real line.
For example, if $\mu$ is supported on the real line,
then Pc$(S(\mu))=S(\mu)$, and if $K$ is a closed interval disjoint
from the support, then there is at most one zero in $K$.
It was shown in Denison--Simon \cite{denissimon} that  if
$x_0\in \R$ is not in the support, then  for some $\d>0$ and all $n$ either
$p_n$ or $p_{n+1}$ has no zero in $(x_0-\d,x_0+\d)$. Note that
if $\mu$ is a symmetric measure on $[-1,-1/2]\cup[1/2,1]$, then
$p_{2n+1}(0)=0$ for all $n$, so the result is sharp.

Any isolated 
point in the support that lies on the outer
boundary attracts precisely one zero. 
Let $z_0$ be an isolated point of $S(\mu)$,
such that its distance from the
convex hull of $S(\mu)\setm\{z_0\}$ is $\d>0$.
Then $p_n$ has at most one zero in the disk $\{|z-z_0|<\d/3\}$
(Simon \cite[Section 8.1]{simon1}).
It is also clear that for any symmetric measure $\mu$ with
$S(\mu)=[-1,-1/2]\cup \{0\}\cup [1/2,1]$ the polynomials
$p_{2n}(\mu)$ have 2 zeros near $0$, so the result is sharp 
(in this case $\d=0$).
Moreover, if $\mu$ lies on the unit circle, then
 there exist two positive constants $C$
and $a$ and a zero $z_n$ of $p_n$ such that $|z_n-z_0|\leq
Ce^{-an}$.

In general, each component of the polynomial convex hull 
consisting of 
more than one point 
attracts infinitely many zeros: if 
$\g$ is a Jordan curve in $\O$ such that 
$S(\mu)\cap \g$ is infinite, then the number of
zeros of $p_n$ that lie inside $\g$ tends to infinity 
(Saff--Totik \cite{safftotik2}).
Mass points of $\mu$ do not necessarily attract zeros
(above we have mentioned that they do if they lie
on the outer boundary).
In fact, it was shown in Saff--Totik \cite{safftotik2} that 
if $\rho$ is the measure on the unit circle given by
the density function $\sin^2(\t/2)$, then for any
measure $\s$ that is supported in the open unit
disk there is a $\l>0$ such that all zeros of 
the $n$-th orthogonal polynomials with respect to $\mu=\rho+\l\s$
tend to the unit circle as $n\to\i$.

Next put a unit mass at every zero of $p_n(\mu)$
(counting multiplicity). This gives
the so called \dword{counting measure} $\nu_{p_n(\mu)}$
on the zero set. \dword{Zero distribution} amounts to finding
the limit behavior of $\frac{1}{n}\nu_{p_n(\mu)}$.
The normalized arc measure on the unit circle (for which $p_n(\mu,z)=z^n$)
shows that if the interior of the polynomial convex hull Pc$(S(\mu))$ is
not empty, then the zeros may be far away from the outer boundary
$\partial \O$, where the equilibrium measure $\o_{S(\mu)}$ is supported.
Thus, assume that  Pc$(S(\mu))$ has empty interior and also that
 there is no Borel set of capacity zero
and full $\mu$-measure, i.e.,
the minimal carrier capacity $c_\mu$  is positive (the $c_\mu=0$ 
case is rather pathological, almost anything can happen with the zeros
then). In this case
\be  \lim_{n\to\i} \k_n^{1/n}=\log\frac1{\c(S(\mu))} \label{limkn}\ee
if and only if
\[\lim \frac{1}{n}\nu_{p_n(\mu)}=\o_{S(\mu)}\]
in weak$^*$ sense,
i.e., asymptotically minimal behavior of $\k_n^{1/n}$ (see
(\ref{lower})) is equivalent to the fact that the zero distribution is
the equilibrium distribution. In a similar way, asymptotic maximal
behavior (see (\ref{cmu})), i.e.,
\be \lim_{n\to\i} \k_n^{1/n}=\frac1{c_\mu}\ee
holds precisely when
\[\lim_{n\to\i} \frac{1}{n}\nu_{p_n(\mu)}=\o_\mu,\]
where $\o_\mu$ is the so called minimal carrier equilibrium measure,
for which a representation like (\ref{greenrep}) is true,
but for the minimal carrier Green function $g_\mu$.

\subsubsection{Regularity}
(\ref{limkn}) is called \dword{regular limit behavior}, and in this case
we write $\mu\in{\bf Reg}$. Thus, the important class
${\bf Reg}$ is defined by the property (\ref{limkn}).
$\mu\in{\bf Reg}$ is equivalent to either of
\begin{itemize}
  \item $\lim_{n\to\i}|p_n(\mu,z)|^{1/n}=\exp({g_{\C\setm S(\mu)}(z,\i)})$,$\ \ \ $$z\not \in {\rm Con}(S(\mu))$
  \item $\limsup_{n\to\i} |p_n(\mu,z)|^{1/n}=1$ for quasi-every $z\in \partial \O$.
\end{itemize}

If $\O$ is a regular set with respect to the Dirichlet problem, then $\mu\in{\bf Reg}$
is equivalent to either of
\begin{itemize}
  \item $\lim_{n\to\i}\|p_n(\mu)\|_{{\rm sup},S(\mu)}^{1/n}=1$
  \item For any sequence $\{P_n\}$ of polynomials of degree $n=1,2,\ldots$
  \[\lim_{n\to\i}\left(\frac{\|P_n\|_{
  {\rm sup}, S(\mu)}}{\|P_n\|_{L^2(\mu)}}\right)^{1/n}=1.\]
\end{itemize}
The last statement expresses the fact that in the $n$-th root sense
the $L^2(\mu)$ and $L^\i$ norms (on $S(\mu)$) are asymptotically
the same.

All equivalent formulations of $\mu\in{\bf Reg}$ point to a
certain ``thickness'' of $\mu$ on its support. Regularity is an
important property, and it is desirable to know ``thickness''
conditions under which it holds. Several regularity criteria are
known, e.g. either of the conditions
\begin{itemize}
  \item all Borel sets $B\subseteq S(\mu)$
  with full measure (i.e with $\mu(B)=\mu(S(\mu))$)
  have capacity $\c(B)=\c(S(\mu))$, i.e., $c_\mu=\c(S(\mu))$ or
  \item $d\mu/d\o_{S(\mu)}>0$ (Radon-Nikodym derivative) $\o_{S(\mu)}$-almost everywhere
\end{itemize}
is sufficient for $\mu\in {\bf Reg}$. Regularity holds under fairly weak assumptions
on the measure, e.g. if $S(\mu)=[0,1]$, and
\[\liminf_{r\to 0}r\log\mu([x-r,x+r])\ge 0\]
for almost every $x\in [0,1]$
(i.e., if $\mu$ is not exponentially small around almost every point), then
$\mu\in{\bf Reg}$.

No necessary and sufficient condition for regularity in terms of the
size of the measure $\mu$ is known.
The only existing necessary condition is for the case $S(\mu)=[0,1]$,
and it reads that for every $\eta>0$
\[\lim_{n\to\i}\c\left(\left\{x\sep\mu([x-1/n,x+1/n])\ge e^{-\eta n}\right\}\right)=\frac14\]
(here 1/4 is the capacity of $[0,1]$).

\sect{Strong,  ratio and weak asymptotics}\label{secreal}
\subsubsection{Strong asymptotics}
Let $\mu$ be supported on $[-1,1]$ and suppose that
the so called \dword{Szeg\H{o} condition}
\be \int_{-1}^1\frac{\log \mu'(t)}{\sqrt{1-t^2}}dt>-\i\label{Szego}\ee
holds, where $\mu'$ is the Radon-Nikodym derivative
of $\mu$ with respect to linear Lebesgue measure.
Note that this condition means that the integral is finite, for it cannot
be $\i$. It expresses a certain denseness of $\mu$, and under this
condition G. Szeg\H{o} proved several asymptotics for the corresponding
orthonormal polynomials $p_n(\mu)$. This theory was developed
on the unit circle and then was translated into the real line.
The \dword{Szeg\H{o} function} associated with $\mu$ is
\be D_\mu(z):=\exp\left(\sqrt{z^2-1}\frac{1}{2\pi}\int_{-
1}^1\frac{\log \mu'(t)}{z-t}\frac{dt}{\sqrt{1-t^2}}\right)
\label{szegofunct}\ee
and it is the outer function in the Hardy space on
$\C\setm [-1,1]$ with boundary values $|D_\mu(x)|^2=\mu'(x)$.
Outside $[-1,1]$ the asymptotic formula
\be p_n(\mu,z)=(1+o(1))\frac{1}{\sqrt{2\pi}}(z+\sqrt{z^2-1})^nD_\mu(z)^{-1}
\label{szegomain}\ee
holds locally uniformly. In particular, the leading coefficient $\k_n$ of
$p_n(\mu)$ is of the form
\be \k_n=(1+o(1))\frac{2^n}{\sqrt{2\pi}}\exp\left(\frac{-1}{2\pi}\int_{-
1}^1\frac{\log \mu'(t)}{\sqrt{1-t^2}}dt\right).\label{szegolead}\ee
If $d\mu(x)=w(x)dx$ and $h(t)=w(\cos t)\sin t$ satisfies a Dini-Lipshitz condition
\[|h(t+\d)-h(t)|\le \frac{C}{|\log \d|^{1+\e}},\qquad \e>0,\]
then with 
\[\G_w(x):=\frac{1}{2\pi}\int_{-1}^1\frac{\log w(\xi)-\log
w(x)}{\xi-x}\left(\frac{1-x^2}{1-\xi^2}\right)^{1/2}d\xi,\]
we have uniformly on $[-1,1]$
\bean (1-x^2)^{1/4}w(x)^{1/2}p_n(x)&=&\left(\frac{2}{\pi}\right)^{1/2}\cos\left(
(n+\frac12)\arccos x+\G_w(x)-\frac{\pi}{4}\right)\\
&+&O((\log n)^{-\e}).\eean
For all these results see Szeg\H{o} \cite{Szego}, Chapter 6.
The Szeg\H{o} condition is also necessary for
these results, e.g. an asymptotic formula
like (\ref{szegomain}) and (\ref{szegolead}) is equivalent
to (\ref{Szego}).

\subsubsection{Ratio asymptotics}
If one assumes weaker conditions then necessarily weaker
results will follow. A large and important class of measures is
 the \dword{Nevai class}
$M(b,a)$ (see Nevai \cite{Nevai}),
for which the coefficients in the three-term recurrence
\[xp_n(x)=a_np_{n+1}(x)+b_np_n(x)+a_{n-1}p_{n-1}(x)\]
satisfy $a_n\to a$, $b_n\to b$. This is equivalent to ratio asymptotics
\[\lim_{n\to\i}\frac{p_{n+1}(z)}{p_n(z)}=\frac{z-b+\sqrt{(z-b)^2-4a^2}}2\]
for large $z$ (actually, away from the support of $\mu$), and the monograph
Nevai \cite{Nevai} contains a very detailed treatment of orthogonal polynomials
in this class. It is also true
that  if the limit of
$p_{n+1}(z)/p_n(z)$ exists at a single non-real $z$, then $\mu\in M(b,a)$
for some $a,b$ (Simon \cite{simon6}).

The classes $M(b,a)$ are scaled versions of each other,
and the most important condition ensuring $M(0,1/2)$ is
given in Rakhmanov's theorem from \cite{rakhmanov}:
if $\mu$ is supported in $[-1,1]$ and
$\mu'>0$ almost everywhere on $[-1,1]$, then $\mu\in M(0,1/2)$. Conversely,
Blumenthal's theorem from \cite{blumenthal}
states that $\mu\in M(0,1/2)$ implies that
the support of $\mu$ is $[-1,1]$ plus at most countably many points
that converge to $\pm1$. Thus, in this respect the extension
of Rakhmanov's theorem given in \cite{denisov} by Denisov
is of importance:
if $\mu'>0$ almost everywhere
on $[-1,1]$ and outside $[-1,1]$ the measure $\mu$ has at most
countably many mass points converging to $\pm1$, then $\mu\in M(0,1/2)$.
However, $M(0,1/2)$ contains many other measures not just those
that are in these theorems, e.g. in
Delyon--Simon--Souillard \cite{simondely} a continuous singular measure 
in the Nevai class
was exhibited, and
the result in Totik \cite{totik1} shows that the Nevai class contains
practically all types of measures allowed by Blumenthal's theorem.

\subsubsection{Weak and relative asymptotics}
Under Rahmanov's condition supp$(\mu)=[-1,1]$, $\mu'>0$ a.e., some
parts of Szeg\H{o}'s theory can be proven in a weaker form (see
e.g. M\'at\'e--Nevai--Totik \cite{mnt1,mnt2}). In these the \dword{Tur\'an determinants}
\[T_n(x):=p_n^2(x)-p_{n-1}(x)p_{n+1}(x)\]
play a significant role. In fact, then given any interval $\D\subset(-1,1)$
the Tur\'an determinant $T_n$ is positive on $\D$ for all large $n$,
and $T_n(x)^{-1}dx$ converges in the weak$^*$ sense to $d\mu$ on $\D$.
Furthermore, the absolutely continuous part $\mu'$ can be
also separately recovered
from $T_n$:
\[\lim_{n\to\i}\int \left|T_n(x)\mu'(x)-\frac{2}{\pi}(1-x^2)^{1/2}\right|dx=0.\]
Under Rahmanov's condition we also
have weak convergence, for example,
\be \lim_{n\to\i}\int f(x)p_n^2(x)\mu'(x)dx =\frac{1}{\pi}
\int_{-1}^1 \frac{f(x)}{\sqrt{1-x^2}}dx\label{rahmweak}\ee
for any continuous function $f$. Pointwise we only know
a highly oscillatory behavior: for almost all $x\in [-1,1]$
\[\limsup_{n\to\i} p_n(x)\ge \frac{2}{\pi}(\mu'(x))^{-1/2}(1-x^2)^{-1/4},\]
\[\liminf_{n\to\i} p_n(x)\le -\frac{2}{\pi}(\mu'(x))^{-1/2}(1-x^2)^{-1/4},\]
and if $E_n(\e)$ is the set of points $x\in [-1,1]$ where
\[|p_n(x)|\ge (1+\e)\frac{2}{\pi}(\mu'(x))^{-1/2}(1-x^2)^{-1/4},\]
then $|E_n(\e)|\to 0$ for all $\e>0$. However, it is not true
that the sequence $\{p_n(\mu,x)\}$ is pointwise
bounded, since
for every $\e>0$
there is a weight function $w>1$ on $[-1,1]$ such that $p_n(0)/n^{1/2-\e}$
is unbounded (see Rakhmanov \cite{rahmanov1}).

Simon \cite{simon6} extended (\ref{rahmweak}) by showing that if
the recurrence coefficients satisfy $b_n\to b$, $a_{2n+1}\to a'$
and $a_{2n}\to a''$, then there is an explicitly calculated
  measure $\r$ depending only
on $b,a',a''$ such that
\be \lim_{n\to\i}\int f(x)p_n^2(x)\mu'(x)dx =
\int_{-1}^1 f(x)d\r(x)\label{simonlim}\ee
for any continuous function $f$, and conversely, if (\ref{simonlim})
exists for $f(x)=x,x^2,x^4$, then $b_n\to b$, $a_{2n+1}\to a'$
and $a_{2n}\to a''$ with some $b,a',a''$.

   For measures in Nevai's class, part of Szeg\H{o}'s theory
can be extended to \dword{relative asymptotics}, i.e., when sequences of orthogonal
polynomials corresponding to two measures are compared. Here is a sample
theorem: let $\a$ be supported in $[-1,1]$ and
in Nevai's class $M(0,1/2)$, and let 
$d\b=gd\a$, where $g$ is a function such that for some polynomial
$R$ both 
$Rg$ and $R/g$ are Riemann integrable. Then
\[\lim_{n\to\i}\frac{p_n(\b,z)}{p_n(\a,z)}=D_g(z)^{-1}\]
uniformly on $\C$ away from $[-1,1]$, where $D_g$ is Szeg\H{o}'s 
function with respect to the measure $g(x)dx$.

\subsubsection{Widom's theory}
Szeg\H{o}'s theory can be extended to measures lying on a single
Jordan curve or arc $J$ (see Kaliaguine \cite{kalyagin}
where also additional outside lying mass points are allowed),
in which case the role of $z+\sqrt{z^2-1}$
in (\ref{szegomain}) is played by the conformal map $\F$
of $\C\setm J$ onto
the exterior of the unit disk, and the role of $2^n$ in
(\ref{szegolead}) is played by the reciprocal of
the logarithmic capacity of $J$
(see Section \ref{secgeneralop}). Things change considerably if the
measure is supported on a set $J$ consisting of
two or more smooth curve or arc components $J_1,\ldots,J_m$.
A general feature
of this case is that $\k_n\c(J)^n$ does not have a limit, its
limit points  fill a whole interval (though if some associated
harmonic measures are all rational then the limit points
may form a finite set). The polynomials themselves
have asymptotic form
\[\frac{p_n(z)}{\k_n}=\c(J)^n\F(z)^n(F_n(z)+o(1))\]
uniformly away from $J$, where $\F$ is the (multi-valued) complex
Green function of the complement $\C\setm J$, and where $F_n$ is
the solution of an $L^2$-extremal problem involving analytic
functions belonging to some class $\G_n$. The functions $F$ in
$\G_n$ are determined by an $H^2$ condition plus an argument
condition, namely if the change of the argument of $\F$ as we go
around $J_k$ is $\g_k 2\pi$ modulo $2\pi$, then in $\G_n$ we
consider functions whose change of the argument around
$J_k$ is $-n\g_k 2\pi$ modulo $2\pi$. Now the point is that these
function classes $\G_n$ change with $n$, and hence so does $F_n$,
and this is the reason that a single asymptotic formula like
(\ref{szegolead}) or (\ref{szegomain}) does not hold. The
fundamentals of the theory were laid out in H. Widom's paper
\cite{widom}; and since then many results have been obtained by F.
Peherstorfer and his collaborators, as well as A. I. Aptekarev, J.
Geronimo, S. P. Suetin and W. Van Assche. The theory has deep connections with
function theory, the theory of Abelian integrals and the theory of
elliptic functions. We refer the reader to the papers
Aptekarev \cite{aptekarev1}, 
Geronimo--Van Assche \cite{geronimo}, 
Peherstorfer \cite{peher2}--\cite{peher4} and
Suetin \cite{Suetin1}--\cite{Suetin2}.

\subsubsection{Asymptotics for Christoffel functions}
The Christoffel functions
\[\l_n(\mu,x)^{-1}=\sum_{k=0}^n p_k(\mu,x)^2\]
behave somewhat more regularly than the orthogonal polynomials.
In M\'at\'e--Nevai--Totik \cite{mnt3} 
it was shown that if $\mu$ is supported on $[-1,1]$,
it belongs to the ${\bf Reg}$ class there (see Section \ref{secgeneralop})
and $\log\mu'$ is integrable over an interval $I\subset [-1,1]$, then for
almost all $x\in I$
\[ \lim_{n\to\i}n\l_n(\mu,x)=\pi\sqrt{1-x^2}\mu'(x).\]
This result is true (see Totik \cite{totik2}) in the form
\[ \lim_{n\to\i}n\l_n(\mu,x)=\frac{d\mu(x)}{d\o_{{\rm supp}(\mu)}(x)},\qquad {\rm a.e.}\ x\in I\]
when the support is a general compact subset of $\R$, $\mu\in {\bf Reg}$
and $\log\mu'\in L^1(I)$.

Often only a rough estimate is needed for Christoffel functions, and such 
a one is
provided in Mastroianni--Totik \cite{mastto}: if $\mu$ is supported on $[-1,1]$ and it is
a doubling measure, i.e.,
\[\mu(2I)\le L\mu(I)\]
for all $I\subset [-1,1]$, where $2I$ is the twice enlarged $I$, then
uniformly on $[-1,1]$
\[\l_n(\mu,x)\sim \mu\left(\Delta_n(x)\right);\qquad \Delta_n(x)=\frac{\sqrt{1-x^2}}{n}+\frac{1}{n^2}.\]

\sect{Recurrence coefficients and spectral measures}
Let $\mu$ be a unit measure of compact 
support on the real line, and
\[ xp_n(x)=a_np_{n+1}(x)+b_np_n(x)+a_{n-1}p_{n-1}(x),\]
the recurrence relation for the corresponding orthogonal
polynomials. We have already mentioned in Section \ref{secwhere}
that $\mu$ is the spectral measure for the Jacobi matrix
\[J= \left(\begin{array}{ccccc}
  b_0 & a_0 & 0 & 0 &\cdots  \\
  a_0 & b_1 & a_1 & 0 & \cdots \\
  0 & a_1 & b_2 & a_2 & \cdots \\
  0 & 0 & a_2 & b_2 &\cdots \\
  \vdots&\vdots&\vdots&\vdots&\ddots
\end{array}\right),\]
and this gives a one-to-one correspondence
between unit measures with compact (and infinite) support
on the real line and Jacobi operators with bounded entries.
Every such Jacobi operator is a bounded self-adjoint operator
on $l_2$, hence operator theory and orthogonal 
polynomials meet at this point, and techniques and questions
from both areas are relevant. Information between the measure
$\mu$ and the recurrence coefficients might be called
the spectral analysis of orthogonal polynomials. In this, 
a special 
role is played by the Chebyshev case $a_n=1/2$ and $b_n=0$, 
where
 the Jacobi matrix is denoted by $J_0$. If the sequences
$\{a_n\}$ and $\{b_n\}$ have limits then we may assume 
$a_n\to1/2$, $b_n\to 0$ (this is just rescaling, and this
was the Nevai class $M(0,1/2)$), hence
in this case the Jacobi operator is a compact
perturbation of $J_0$,
and one of the main questions of the theory is
how properties of $J-J_0$ are reflected in the 
spectral measure $\mu$.

We have already mentioned in Section \ref{secreal} that
$\mu\in M(0,1/2)$ implies that the support of $\mu$ is 
$[-1,1]$ plus some additional mass points converging to $\pm1$
(called Blumenthal's theorem in orthogonal polynomials;
it is a special case of  Weyl's theorem in operator theory on the
invariance of the essential spectrum under
compact perturbation).
Conversely, if the support of $\mu$ is 
$[-1,1]$ plus some additional mass points converging to $\pm1$ 
and $\mu'(x)>0$ for almost all $x\in [-1,1]$, then
$\mu\in M(a,b)$ (Denisov \cite{denisov}, Rakhmanov \cite{rakhmanov},
Nevai--Totik \cite{nevtot1}). 
No spectral characterization of $\mu\in M(0,1/2)$ is known;
this important class seems to contain all sorts of measures.
For example, if $\nu$ is any measure
with support $[-1,1]$ then there is a $\mu\in M(0,1/2)$
which is absolutely continuous with respect to $\nu$. In particular,
$M(0,1/2)$ contains discrete measures, continuously singular measures
or measures that are given by a continuous density which is positive
on a set of measure $<\e$.

Strengthening the condition $\mu\in M(0,1/2)$ can
be done in several ways. After numerous works in the
subject 
by Szeg\H{o}, Shohat, Geronimus, Krein, Kolmogorov and others,
a complete
characterization for $J-J_0$ being a Hilbert-Schmidt operator
was given in \cite{killip} 
by R. Killip and B. Simon (note that $\mu$ is assumed to have
total mass 1):
\be \sum_n(a_n-1/2)^2+\sum_n b_n^2<\i\label{square}\ee
if and only if the following conditions hold:
\begin{description}
\item[(i)] the support of $\mu$ is 
$[-1,1]$ plus some additional mass points $E^\pm_j$ converging to $\pm1$,
\item[(ii)] if $\mu'$ is the absolutely continuous part of $\mu$
on $[-1,1]$, then
\[\int_{-1}^1(\log\mu'(t))\sqrt{1-t^2}dt>-\i,\]
\item[(iii)] for the mass points $E^\pm_j$ lying outside $[-1,1]$ we have
\[\sum_j|E^+_j-1|^{3/2}+\sum_j|E^-_j+1|^{3/2}<\i.\]
\end{description}

It was also shown in Killip--Simon \cite{killip} that if
$J-J_0$ is trace class, i.e., 
\[\sum_n|a_n-1/2|+\sum_n b_n<\i,\]
then Szeg\H{o}'s condition
\be \int_{-1}^1 \frac{\log \mu'(t)}{\sqrt{1-t^2}}dt>-\i
\label{szego3}\ee
holds. The conclusion is also true if $\mu\in M(0,1/2)$,
for the mass points $E^\pm_j$ lying outside $[-1,1]$ we have
\[\sum_j|E^+_j-1|^{1/2}+\sum_j|E^-_j+1|^{1/2}<\i,\]
and 
\be \limsup_n(2^na_1\cdots a_n)>0.\label{lsup}\ee
If the support of $\mu$ is contained in $[-1,1]$, then
Szeg\H{o}'s condition automatically holds if (\ref{lsup}) is true.
Actually, when supp$(\mu)=[-1,1]$ then
 Szeg\H{o}'s condition
(\ref{szego3}) is equivalent to (\ref{square}) and to 
the (conditional) convergence
of the series $\sum_n(a_n-1/2)$ and $\sum_nb_n$.

There is also an extended theory of orthogonal polynomials
with several different applications when the recurrence coefficients
do not converge, but they are asymptotically periodic in the
sense that for some $k$ all the sequences $(a_{kn+j})_{n=1}^\i$
and $(b_{kn+j})_{n=1}^\i$, $j=1,\ldots,k$ converge. These are related
to so called sieved orthogonal polynomials and to orthogonal polynomials
generated by polynomial mappings. In this case
the essential support of the spectral measure
lies on several intervals. There are numerous
papers on this subject by M. E. H. Ismail, N. A. Al-Salam,
J. A. Charris, J. Wimp, J. Bustoz, J. Geronimo, W. Van Assche,
F. Peherstorfer, R. Steinbauer,  N. I. Akhiezer, B. P. Osilenker and others;
see e.g. Charris--Ismail \cite{ismail5}, 
Geronimo--Van Assche \cite{geronimo},
Peherstorfer \cite{peher3}, Peherstorfer--Steinbauer \cite{peher1}, 
Akhiezer \cite{akhiezer1} 
for details and for further references.

\sect{Exponential and Freud weights}\label{secfreud} These are
weight functions of the form $e^{-2Q(x)}$, where $x$ is on the
real line or on some subinterval thereof. For simplicity we shall
first assume that $Q$ is even. We get \dword{Freud weights} when
$Q(x)=|x|^\a$, $\a>0$, $x\in \R$, and \dword{Erd\H{o}s weights} if $Q$
tends to infinity faster than any polynomial as $|x|\to\i$. G.
Freud started to investigate these weights in the sixties and
seventies, but they independently appeared also in the Russian
literature and in statistical physics. One can safely say that
some of Freud's problems and the work of P. Nevai and E. 
A. Rahmanov
were the primary cause of the sudden revitalization of the theory
of orthogonal polynomials since the early 1980's.  In the last 20
years D. Lubinsky with coauthors have conducted systematic studies
on exponential weights, see e.g. Levin--Lubinsky \cite{Levlub, Levlub1}, 
Lubinsky \cite{Lub},
Lubinsky--Saff \cite{Lubsaff}, Van Assche \cite{Assche}; 
we should mention the names
E. Levin, E. B. Saff, W. Van Assche, E. A. Rahmanov and
H. N. Mhaskar. In the mid 1990's a new stimulus came from the
Riemann--Hilbert approach that was used together with the steepest
descent method by P. Deift, T. Kriecherbauer, 
K. T.-R. McLaughlin, S. Venakides
and X. Zhou (\cite{Deift1})
to give complete asymptotics when $Q$ is analytic.

One can roughly say that because of the fast vanishing of the
weight around infinity, things happen on a finite subinterval $[-a_n,a_n]$
(depending on the degree of the polynomials), and on $[-a_n,a_n]$
techniques developed for $[-1,1]$ are applied. For Freud weights one
can also make the substitution $x\to n^{1/\l} x$ and go to
orthogonality with respect to the varying weight $e^{-n|x|^\l}$,
in which case things are automatically reduced to a finite interval
which is the support of a weighted energy problem.

The $a_n$ are the so called \dword{Mhaskar-Rahmanov-Saff numbers}
defined by
\be n=\frac{2}{\pi}\int_0^1\frac{a_ntQ'(a_nt)}{\sqrt{1-
t^2}}\;dt.\label{mrs}\ee
The zeros of $p_n(w^2)$, $w(x)=\exp(-Q(x))$  are spreading out
and the largest zero is very close to $a_n$, which tends to $\i$.

To describe the distribution of the zeros and the behavior of the polynomials
one has to make appropriate contractions. Let us consider first
the case of Freud weight $w(x)=\exp (-|x|^\a)$, and let $p_n$ be
the $n$-th orthogonal polynomial with respect to $w^2$ (on $(-\i,\i))$.
In this case
\[ a_n=n^{1/\a}\g_\a, \qquad \gamma _\alpha  ^{\a}:=
\Gamma \left({\alpha \over 2}\right)\Gamma \left({1\over
2}\right)\Bigl/2\Gamma\left({\alpha \over 2} +
{1\over 2}\right).  \]
Thus, for the largest zero $x_{n,n}$ we have
$x_{n,n}/n^{1/\a}\to \g_\a$ as $n\to\i$, and to describe
zero distribution we divide (contract) all zeros $x_{n,i}$ by
$n^{1/\a}\g_\a$. These contracted zeros asymptotically
have the \dword{Ullman distribution}
\be \frac{d\mu_w(t)}{dt}
:= {\alpha \over \pi } \int  ^1_{|t|} \frac{u^{\alpha -1}}{\sqrt{u^2 -t^2}} \, du
,\quad t \in  [-1,1].\label{ullm}\ee
This measure $\mu_w$ minimizes the weighted energy
\be \int\int\log\frac{1}{|x-t|}d\mu(x)d\mu(t)+2\int Q d\mu
\label{equiw}\ee
among all probability measures compactly supported on $\R$.
It is a general feature of exponential weights that
the behavior of zeros of the polynomials is governed by
the solution of a weighted energy problem (weighted equilibrium
measures, see Saff--Totik \cite{safftot}). If $\k_n$ is the leading coefficient of
$p_n$, i.e., $p_n(z)=\k_nz^n+\cdots$, then (Lubinsky--Saff \cite{Lubsaff})
\[\lim_{n\to\infty} \k_n\pi^{1/2}2^{-
n}e^{-n/\a}n^{(n+1/2)/\a}= 1,\]
and we have
\bean &&\lim_{n\to\i}|p_n(n^{1/\alpha }\gamma_\alpha z)|^{1/n}\\
&&= \exp \left(  \log |z + \sqrt{z^2-1}|
+ {\rm Re} \int
^1_0
{zu^{\alpha -1}\over \sqrt{z^2-u^2}}\, du \right)\eean
locally uniformly outside $[-1,1]$. This latter 
is so called $n$-th root asymptotics, while the
former is strong asymptotics. Strong asymptotics for $p_n(z)$
on different parts of the
complex plane was given using the Riemann--Hilbert approach,
see Deift \cite{Deift1} and 
Kriecherbauer--McLaughlin \cite{Kriecher} and the references there.
On the real line we have a Plancherel--Rotach type formula
\[ n^{1/2\a}p_n(w_\a;n^{1/\a}\g_\a x)\exp(-n\g_\a^\a|x|^\a)-\]
\[-\sqrt\frac{2}{\pi}
\frac{1}{\sqrt[4]{1-x^2}}\cos\left(\frac{1}{2}\arccos x+n\pi
\mu_w([x,1])
-\frac{\pi}{4}\right)\to 0\]
uniformly on any subinterval of $(-1,1)$.

Things become more complicated for non-Freud weights, but the
corresponding results are of the same flavor. In this case
the weight is not necessarily symmetric, but under some
conditions (like $Q$ being convex or $xQ'(x)$ being increasing
 for
$x>0$ and an analogous condition for $x<0$) the relevant
weighted equilibrium measure's support is an interval, and
the definition of the
Mhaskar--Rahmanov--Saff numbers $a_{\pm n}$ is
\bean n&=&\frac{1}{\pi}\int_{a_{-n}}^{a_n}\frac{xQ'(x)}{\sqrt{(x-a_{-n})(a_n-x)}}dx,\\
0&=&\frac{1}{\pi}\int_{a_{-n}}^{a_n}\frac{Q'(x)}{\sqrt{(x-a_{-n})(a_n-x)}}dx.\eean
Now one solves the weighted equilibrium problem (\ref{equiw})
for all measures $\mu$ {\it with total mass $n$}, and if $\mu_n$ is the solution
then $[a_{-n},a_n]$ is the support of $\mu_n$ and $\mu_n/n$ will
play the role of the measure $\mu_w$ from (\ref{ullm}) above.

The weight does not even have to be defined on all
$\R$, e.g. in \cite{Levlub} a theory was developed 
by Levin and Lubinsky that
simultaneously includes far reaching generalizations of
non-symmetric Freud, Erd\H{o}s and Pollaczek weights
 such as

(a) nonsymmetric Freud-type weights
\[Q(x)=\left\{\begin{array}{cc}
  |x|^\a, & x\in[0,\i) \\
  |x|^\b, & x\in(-\i,0), 
\end{array}\right.\]

(b) nonsymmetric Erd\H{o}s weights such as 
\[Q(x)=\left\{\begin{array}{cc}
  \exp_l(|x|^\a)-\exp_l(0), & x\in[0,\i) \\
  \exp_k(|x|^\b)-\exp_k(0), & x\in(-\i,0) 
\end{array}\right.\]
with $\exp_l$ is the $l$-times iterated exponential function, or

(c) nonsymmetric Pollaczek type weights that vanish fast at $\pm1$ such as 
\[Q(x)=\left\{\begin{array}{cc}
  \exp_l((1-x)^{-\a})-\exp_l(1), & x\in[0,1) \\
  \exp_k((1-x)^{-\b})-\exp_k(0), & x\in(-1,0). 
\end{array}\right.\]

In all cases the interval $[a_{-n},a_n]$ is where things happen, 
e.g. this is the shortest
interval on which the supremum norm of a weighted polynomial is attained:
\[\|wP_n\|_{\sup}=\|wP_n\|_{\sup,[-a_{-n},a_n]}\]
for all polynomials of degree at most $n$. These numbers $a_{\pm n}$
are everywhere in
the theory, e.g. 
\[\sup_x |p_n(x)|w(x)|x-a_n|^{1/4}|x-a_{-n}|^{1/4}\sim 1.\]

\sect{Sobolev orthogonality}\label{secsobolev} In 
\dword{Sobolev orthogonality} we consider orthogonality with respect to an inner
product \be (f,g)=\sum_{k=0}^r\int
f^{(k)}\overline{g^{(k)}}d\mu_k, \label{sob}\ee where $\mu_k$ are
given positive measures. There are several motivations for this
kind of orthogonality. Perhaps the most natural one is smooth data
fitting. The Spanish school around F. Marcell\'an, G. Lopez and A.
Martinez-Finkelshtein has been particularly active in developing
this area (see the surveys Marcell\'an--Alfaro--Rezola \cite{Marcellan} 
and Martinez-Finkelshtein \cite{Martinez,
Martinez1} and the references therein).

In this section let $Q_n(z)=z^n+\cdots$ denote the {\it monic} orthogonal
polynomial with respect to the Sobolev inner product (\ref{sob}), and
$q_n(\mu_k)$ the monic orthogonal polynomials with respect to the measure
$\mu_k$.

Most arguments for the standard theory fail in this case, e.g. it is no
longer true that the zeros lie in the convex hull of the support of the
measures $\mu_k$, $k=0,1,\ldots,r$. It is not even known if
the zeros are bounded if all the measures $\mu_k$ have compact support.
Nonetheless, for the case $r=1$, and $\mu_0,\mu_1\in{\bf Reg}$
(see Section \ref{secgeneralop}) it was shown in 
Gautschi--Kuijlaars \cite{Kuijl} that
the asymptotic distribution of the zeros of the {\it derivative} $Q_n'$
is the equilibrium measure $\o_{E_0\cup E_1}$,
where $E_i$ is the  support of $\mu_i$, $i=0,1$
(which also have to be assumed to be regular). Furthermore, if, in addition,
$E_0\subseteq E_1$, then the asymptotic zero distribution of $Q_n$ is
$\o_{E_0}$.

In general, both the algebraic and the asymptotic/analytic
situation is quite complicated, and there are essentially two important cases
which have been understood to a satisfactory degree.
\medskip

\noindent{\it Case I: The discrete case}. In this case $\mu_0$ is
some ``strong'' measure, e.g. from the Nevai class $M(b,a)$ (see
Section \ref{secreal}), and $\mu_1,\ldots,\mu_k$ are finite
discrete measures. It turns out that then the situation is similar
to adding these discrete measures to $\mu_0$ (the new measure will
also be in the same Nevai class), and considering standard
orthogonality with respect to this new measure. For ecample, if $r=1$,
then
\[\lim_{n\to\i}\frac{Q_n(z)}{q_n(\mu_0+\mu_1,z)}=1\]
holds uniformly on compact subsets of $\C\setm {\rm supp}(\mu_0+\mu_1)$.
Thus, the Sobolev orthogonal polynomials differ from those
of the measure $\mu_0$, but not more than what happens when adding mass points
to $\mu_0$.

In this discrete case the $Q_n$'s satisfy a higher order recurrence
relation, hence this case is also related to matrix orthogonality
(see the end of the Section \ref{secmatrix}).
\medskip

\noindent{\it Case II: The Szeg\H{o} case}. Suppose now that $\mu_0,\ldots,\mu_k$
are all supported on the same smooth curve or arc $J$,
and they satisfy Szeg\H{o}'s condition
there (see Section \ref{secreal}). In this case the $k$-th derivative of $Q_n$ satisfies,
locally uniformly in the complement of $J$, the
asymptotic formula
\[\lim_{n\to\i}\frac{Q_n^{(k)}(z)}{n^kq_{n-k}(\mu_r,z)}=\frac{1}{[\F'(z)]^{m-k}},\]
where $\F$ is the conformal map that maps $\C\setminus J$ onto the complement
of the unit disk. That is, in this case the measures $\mu_0,\ldots,\mu_{r-1}$ do
not appear in the asymptotic formula, only $\mu_r$ matters. The reason for this
is the following: $Q=Q_n$ minimizes
\be (Q,Q)=\sum_{k=0}^r \int |Q^{(k)}|^2d\mu_k\label{mineq}\ee
among all monic polynomials of degree $n$, while
$q=q_{n-k}(\mu_k)$ minimizes
\[\int |q|^2d\mu_k\]
among all monic polynomials of degree $n-k$. But the polynomial
$Q_n^{(k)}(t)=n(n-1)\cdots (n-k+1)t^{n-k}+\cdots$ is a monic polynomial
times the factor $n(n-1)\cdots (n-k+1)\sim n^k$, and this factor
is dominant for $k=r$, so everything else will be negligible. There are results
for compensation of this $n^k$ factor which lead to
Sobolev orthogonality with respect to varying measures.

Under the much less restrictive assumption that $\mu_0\in {\bf Reg}$
(see Section \ref{secgeneralop})
and the other measures $\mu_k$ are supported in the support $E$ of
$\mu_0$ it is true (L\'opez--Pijeira-Cabrera--Izquierdo \cite{Lopez})
that the asymptotic
zero distribution of $Q_n^{(k)}$ is
the equilibrium measure $\o_E$ for all $k$,
\[\lim_{n\to\i}\|Q_n^{(k)}\|_{{\rm sup}, E}^{1/n}=\c(E),\]
and hence, away from the zeros in the unbounded component of the complement of $E$,
we have
\[\lim_{n\to\i}|Q_n^{k}(z)|^{1/n}=e^{g_{\C\setm E}(z)}\]
where $g_{\C\setm E}$ is the Green function for this unbounded component.

The techniques developed for exponential weights and for Sobolev
orthogonality were combined in Geronimo--Lubinsky--Marcellan 
\cite{gerlubmarc} to prove strong
asymptotics for Sobolev orthogonal polynomials when $r=1$ and
$\mu_0=\mu_1$ are  exponential weights.

\sect{Non-Hermitian orthogonality}\label{secnonherm}
We refer to \dword{non-Hermitian orthogonality} in either of these
cases:

\begin{itemize}
  \item the measure $\mu$ is non-positive or even complex-valued
and we consider $p_n$ with
\be \int p_n(z) \overline{z^k}d\mu=0, \qquad k=0,1,\ldots,n-1,\label{ort1}\ee
  \item $\mu$  is again non-positive or complex-valued,
or positive but lies on a complex curve or arc
and orthogonality
is considered without complex conjugation, i.e., 
\be \int p_n(z) z^kd\mu=0, \qquad k=0,1,\ldots,n-1.\label{ort2}\ee
\end{itemize}
More generally, one could consider non-positive inner products,
but we shall restrict our attention to complex measures
and orthogonality (\ref{ort2}).

As an example, consider the diagonal Pad\'e approximant to
the Cauchy transform
\[f(z)=\int\frac{d\mu(t)}{z-t}\]
of a signed or complex-valued measure, i.e., consider polynomials
$p_n$ and $q_n$ of degree at most $n$
such that
\[f(z)p_n(z)-q_n(z)=O(z^{-n-1})\]
at infinity. Then $p_n$ satisfies the non-Hermitian
orthogonality relation
\be \int p_n(x)x^jd\mu(x)=0,\qquad j=0,1,\ldots,n-1.
\label{stahlp}\ee

In this non-Hermitian
case even the Gram-Schmidt 
orthogonalization process may fail, and
then $p_n$ is defined as the solution of the orthogonality condition
(\ref{ort1}), resp. (\ref{ort2}), which give a system of homogeneous
equations for the coefficients of $p_n$. Thus, $p_n$ may have
smaller degree
than $n$, and things can get pretty wild with this
kind of orthogonality. For example, in the simple case
\[d\mu(x)=(x-\cos\pi\a_1)(x-\cos\pi\a_2)(1-x^2)^{-1/2}dx, \quad x\in [-1,1],\]
with $0<\a_1<\a_2<1$ rationally independent algebraic numbers,
the zeros of $p_n$ from (\ref{stahlp}) are dense
on the whole complex plane (compare this
with the fact that for positive $\mu$ all zeros lie in $[-1,1]$).
In Stahl \cite{stahl1} it was shown
that it is possible to construct a complex measure $\mu$ on
$[-1,1]$, such that for an arbitrary prescribed asymptotic
behavior some subsequence $\{p_{n_k}\}$ will have this zero behavior.
Nonetheless, the asymptotic distribution of the zeros
is again the equilibrium distribution of the support of $\mu$
under regularity conditions on $\mu$. For example, this
is the case if 
\begin{itemize}
\item $|\mu|$ belongs to the ${\bf Reg}$ class (see
Section \ref{secgeneralop}),
and the argument of $\mu$, i.e., $d\mu(t)/d|\mu|(t)$, is of bounded variation
(Baratchart--K\"ustner--Totik \cite{Baratchart}), or
  \item $d\mu(x)=g(\arccos x)(1-x^2)^{-1/2}dx$, $x\in [-1,1]$, $g$ is bounded 
  away from zero and infinity, and satisfies $|g(\t+\d)-g(\t)|\le K|\log\d|^{-1-\d}$, or
  \item $\mu$ is supported on finitely many intervals, the argument of $\mu$ is uniformly
  continuous and for $a(\d)=\inf_{x\in{\rm supp}(\mu)}|\mu|([x-\d,x+\d])$ 
  the property $\lim_{\d\to0} \log a(\d)=0$ holds (Stahl \cite{stahl1}).
\end{itemize}
In \cite{stahl3}--\cite{stahl1} H. Stahl obtained
asymptotics for non-Hermitian orthogonal polynomials even for varying
measures and gave several
applications of them to Pad\'e approximation. When
the measure $\mu$ is of the form $d\mu(x)=g(x)(1-x^2)^{-1/2}dx$, $x\in [-1,1]$,
with an analytic $g$, for $z\in \C\setm [-1,1]$,
a strong asymptotic formula
of the form
\[\frac{p_n(z)}{\k_n}=(1+o(1))
\frac{(z+\sqrt{z^2-1})^n}{2^n}D_\mu(z)^{-1}\exp\left(\frac{1}{2\pi}\int_{-
1}^1\frac{\log \mu'(t)}{\sqrt{1-t^2}}dt\right)\]
(with $D_\mu$ the Szeg\H{o} function (\ref{szegofunct}))
was proved by  J. Nuttall \cite{Nuttall1},
\cite{Nuttall2}, A. A. Gonchar and S. P. Suetin \cite{Goncharsuetin}.
For a recent Riemann--Hilbert
approach see the paper \cite{aptekarev2}
by A. I. Aptekarev and W. Van Assche.
A similar result holds on the support
of the measure, as well as for the case of varying weights,
see Aptekarev--Van Assche \cite{aptekarev2}.

\sect{Multiple orthogonality}\label{secmultiple}
\def\und{{\underline n}}
Multiple orthogonality comes from simultaneous Pad\'e approximation.
It is a relatively new area where we have to mention the names of
E. M. Nikishin, V. N. Sorokin, A. A. Gonchar and E. A.
Rahmanov, A. I. Aptekarev, A. B. J. Kuijlaars, J. Geronimo and W. Van Assche
(see the survey \cite{assche2} by W. Van Assche and the references there
and the paper Gonchar--Rakhmanov \cite{Goncharrahmanov}).
The analogues of many classical concepts and properties have
been found, and also  the analogues of the classical orthogonal polynomials are known,
e.g. in the multiple Hermite case
the measures are $d\mu(x)=e^{-x^2+c_jx}dx$.

Asymptotic behavior of multiple orthogonal polynomials is not
fully understood yet due to the interaction of the different
measures. For the existing results see Aptekarev \cite{aptekarev},
Van Assche \cite{assche2},
Van Assche's Chapter 23 in \cite{Ismail}
and the references therein.

\subsubsection{Types and normality}
On $\R$ let there be given $r$ measures $\mu_1,\ldots,\mu_r$ with finite moments
and infinite support,
and consider multiindices $\und =(n_1,\ldots,n_r)$ of nonnegative integers
with norm $|\und |=n_1+\cdots+n_r$. There are two types of multiple
orthogonality corresponding to
the appropriate Hermite-Pad\'e approximation.

In \dword{type I  multiple orthogonality}
 we are looking for polynomials $Q_{\und ,j}$ of degree $n_j-1$
for each $j=1,\ldots,r$, such that
\[\sum_{j=1}^r \int x^k Q_{\und ,j}(x)d\mu_j(x)=0,\qquad k=0,1,\ldots,|\und |-2.\]
These orthogonality relations give $|\und |-1$ homogeneous linear equations for
the $|\und|$ coefficients of the $r$ polynomials $Q_{\und ,j}$, so there is a
non-trivial solution. If the rank of the system is $|\und |-1$, then the solution
is unique up to a multiplicative factor, in which case the index $\und$ is called
\dword{normal}. This happens precisely if each $Q_{\und ,j}$ is of exact degree $n_j-1$.

In \dword{type II multiple orthogonality} we are looking 
for a single polynomial $P_{\und} $ of degree
$|\und|$ such that
\bean \int x^kP_\und (x)d\mu_1(x)&=&0,\qquad k=1,\ldots,n_1-1\\
&\vdots&\\
\int x^kP_\und (x)d\mu_r(x)&=&0,\qquad k=1,\ldots,n_r-1.\eean
These are $|\und|$ homogeneous linear equations for the $|\und|+1$
coefficients of $P_\und$, and again if the solution is unique up
to a multiplicative constant, then $\und$ is called  normal. This is
again equivalent to $P_\und$ being of exact degree $\und$.

$\und$ is normal for type I orthogonality precisely when it
is normal for type II, so we just speak of normality.
This is the case,
for example, if the $\mu_j$'s are supported on intervals $[a_j,b_j]$
that are disjoint except perhaps for their endpoints;
in fact, in this case $P_\und$ has $n_j$ simple zeros on $(a_j,b_j)$.
Normality also holds if $d\mu_j=w_jd\mu$ with a common
$\mu$ supported on some interval $[a,b]$,
and  for all $m_j\le n_j$, $j=1,\ldots,r$,  
every non-trivial linear combination of the functions 
\[w_1(x),xw_1(x),\ldots,x^{m_1-1}w_1(x),w_2(x),xw_2(x),\ldots,x^{m_r-1}w_r(x)\]
has at most $m_1+\cdots +m_r-1$ zeros on $[a,b]$ (this means that these
functions form a so called Chebyshev system there).
In this case $P_\und$ has $|\und|-1$ zeros on $[a,b]$.

\subsubsection{Recurrence formulae}
To describe recurrence formulae, let $\underline e_j=(0,\ldots,1,\ldots,0)$
where the single 1 entry is at position $j$. Under 
the normality assumption
if  $P_\und$ is the monic orthogonal polynomial, then for any $k$
\[xP_\und(x) =P_{\und+\underline e_k}(x)+
a_{\und ,0}P_{\und}(x)+\sum_{j=1}^r
a_{\und ,j}P_{\und-\underline \e_j}(x).\]
Another recurrence formula is
\[xP_\und(x) =P_{\und+\underline e_k}(x)+
b_{\und,0}P_{\und}(x)+\sum_{j=1}^r
b_{\und,j}P_{\und-\underline \e_{\pi(1)}-\cdots-\underline e_{\pi(j)}}(x),\]
where $\pi(1),\ldots,\pi(r)$ is an arbitrary, but fixed, permutation of
$1,2,\ldots,r$.
The orthogonal polynomials with different indices are strongly related
to one another, e.g. $P_{\und+\underline e_k}(x)-P_{\und+\underline e_l}(x)$
is a constant multiple of $P_{\und}(x)$.

If $d\mu_j=w_jd\mu$, then similar recurrence relations hold
in case of type I orthogonality for
\[Q_\und(x)=\sum_{j=1}^r Q_{\und,j}(x)w_j(x).\]

Also, type I and type II are related by a \dword{biorthogonality} property:
\[\int P_\und Q_{\underline m}d\mu=0\]
except for the case when $\underline m=\und+\underline e_k$ for some
$k$, and then the previous integral is not zero (under the
normality condition).

\def\und{{\underline m}}
To describe an analogue of the Christoffel-Darboux formula
let $\{\und_j\}$ be  a sequence of multiindices such that $\und_0$ is
the identically 0 multiindex, and $\und_{j+1}$ coincides with $\und_j$
except for one component which is 1 larger than the corresponding component
of $\und_j$. Set $P_j=P_{\und_j}$, $Q_j=Q_{\und_{j+1}}$ and  with $\und=\und_{n}$
\[h_\und^{(j)}:=\int P_\und(x)x^{(\und)_j}d\mu_j(x),\]
where $(\und)_j$ denotes the $j$-th component of the multiindex $\und$.
Then (see Daems--Kuijlaars \cite{damems}), again with $\und=\und_{n}$,
\[(x-y)\sum_{k=0}^{n-1}P_k(x)Q_k(y)=P_\und(x)Q_\und(y)-\sum_{j=1}^r
\frac{h_\und^{(j)}}{h_{\und-\underline e_j}^{(j)}}P_{\und-\underline e_j}(x)
Q_{\und+\underline e_j}(y).\]
Thus, the left hand side depends only on $\und=\und_n$ and not on the particular
choice of the sequence $\und_j$ leading to it.

\def\und{{\underline n}}
\subsubsection{The Riemann--Hilbert problem}
There is an approach (see Van Assche--Geronimo--Kuijlaars \cite{assche1})
to both types of multiple orthogonality
in terms of matrix-valued Riemann--Hilbert problem for $(r+1)\times (r+1)$ matrices
$Y=(Y_{ij}(z))_{i,j=0}^r$.

If $d\mu_j(x)=w_jdx$, then one requires that 

\begin{itemize}
  \item $Y$ is analytic on $\C\setm \R$,
  \item if $Y^\pm(x)$ denote the limit of $Y(z)$ as $z\to x\in \R$ 
from the upper, respectively the lower, half plane, then 
we have $Y^+(x)=Y^-(x)S(x)$,
where
\[S(x):=\left[\begin{array}{ccccc}
  1& w_1(x) & w_2(x) & \cdots & w_r(x) \\
  0 &1 & 0 & \cdots & 0 \\
  0 & 0 & 1 & \cdots & 0 \\
  \vdots & \vdots  & \vdots  & \ddots & \vdots \\
  0 & 0 & 0 & \cdots & 1 
\end{array}\right],\]
  \item as $z\to\i$
\[Y(z)=\left(I+O\left(\frac{1}{z}\right)\right)
\left[\begin{array}{ccccc}
  z^{|\und|}& 0 & 0 & \cdots & 0 \\
  0 &z^{-n_1} & 0 & \cdots & 0 \\
  0 & 0 & z^{-n_2} & \cdots & 0 \\
  \vdots & \vdots  & \vdots  & \ddots & \vdots \\
  0 & 0 & 0 & \cdots & z^{-n_r} 
\end{array}\right].\]

\end{itemize}
The  first entry $Y_{11}(z)$ is precisely the orthogonal polynomial
$P_\und$ of type II, and the other entries are also explicit in terms
of the $P_\und$'s and $w_j$'s (all other entries
are either a constant multiple of $P_{\und-\underline e_k}$ or a Cauchy transform
of its multiple with $w_j$). For type I orthogonality the transfer
matrix is
\[S(x):=\left[\begin{array}{ccccc}
  1& 0 & 0 & \cdots & 0 \\
  -w_1(x) &1 & 0 & \cdots & 0 \\
  -w_2(x) & 0 & 1 & \cdots & 0 \\
  \vdots & \vdots  & \vdots  & \ddots & \vdots \\
  -w_r(x) & 0 & 0 & \cdots & 1 
\end{array}\right],\]
the behavior at infinity is of the form
\[Y(z)=\left(I+O\left(\frac{1}{z}\right)\right)
\left[\begin{array}{ccccc}
  z^{-|\und|}& 0 & 0 & \cdots & 0 \\
  0 &z^{n_1} & 0 & \cdots & 0 \\
  0 & 0 & z^{n_2} & \cdots & 0 \\
  \vdots & \vdots  & \vdots  & \ddots & \vdots \\
  0 & 0 & 0 & \cdots & z^{n_r} 
\end{array}\right],\]
and the multiple orthogonal polynomials $Q_{\und,j}$
are $Y_{1,j+1}/2\pi i$.

\sect{Matrix orthogonal polynomials}\label{secmatrix}
In the last 20 years the fundamentals of matrix orthogonal polynomials have been
developed mainly by A. Dur\'an and his coauthors
(see also the work \cite{Aptekarevnikishin} by
A. I. Aptekarev and E. M. Nikishin). The theory shows many
similarities with the scalar case, but there
is an unexpected  richness which is still to be explored.

For all the results in this section see 
L\'opez-Rodriguez--Dur\'an \cite{duran} and 
Dur\'an--Gr\"unbaum \cite{duran1}
and the numerous references there.

\subsubsection{Matrix orthogonal polynomials}
An $N\times N$ matrix
\[P(t)=\left(\begin{array}{ccc}
  p_{11}(t) & \cdots & p_{1N}(t) \\
  \vdots & \ddots & \vdots \\
  p_{N1}(t) & \cdots & p_{NN}(t)
\end{array} \right)\]
with polynomial entries $p_{ij}(t)$ of degree at most
$n$ is called a \dword{matrix polynomial} of degree at most $n$.
Alternatively, one can write
\[P(t)=C_nt^n+\cdots+C_0\]
with numerical matrices $C_n,\ldots,C_0$ of size $N\times N$.

The number $t=a$ is called a \dword{zero} of $P$ if $P(a)$ is singular, 
and the \dword{multiplicity}
of $a$ is the multiplicity of $a$ as a zero of det$P(a)$. When the 
\dword{leading coefficient matrix $C_n$}
is non-singular, then $P$ has $nN$ zeros counting multiplicity.

From now on we fix the dimension to be $N$, but the degree $n$ can be
any natural number. $I$ will denote the $N\times N$ unit matrix and $0$ stands
for all kinds of zeros (numerical or matrix).

A \dword{matrix}
\[W(t)=\left(\begin{array}{ccc}
  \mu_{11}(t) & \cdots & \mu_{1N}(t) \\
  \vdots & \ddots & \vdots \\
  \mu_{N1}(t) & \cdots & \mu_{NN}(t)
\end{array} \right)\]
\dword{of complex measures} defined on (or part of) the real line
is positive definite if for any Borel set $E$ the numerical
matrix $W(E)$ is positive semidefinite.
We assume that all moments of $W$ are finite. With such a matrix
we can define a matrix inner product on the space of
$N\times N$ matrix polynomials via
\[(P,Q)=\int P(t)dW(t)Q^*(t),\]
and if $(P,P)$ is nonsingular for any $P$ with
nonsingular leading coefficient, then just
as in the scalar case one can generate a sequence
$\{P_n\}_{n=0}^\i$ of matrix
polynomials of degree $n=0,1,\ldots$
which are orthonormal  with respect to $W$:
\[\int P_n(t)dW(t)P_m^*(t)=\left\{
\begin{array}{ll}
  0 & \mbox{if $n\not=m$} \\
  I & \mbox{if $n=m$,}
\end{array}\right.\]
and here $P_n$ has nonsingular leading coefficient matrix.
The sequence $\{P_n\}$ is determined only up to left multiplication
 by unitary matrices, i.e., if $U_n$ are unitary matrices,
then the polynomials $U_nP_n$ also form an orthonormal system with
respect to $W$.

\subsubsection{Three-term recurrence and quadrature}
Just as in the scalar case, these orthogonal polynomials satisfy
a three-term recurrence relation
\be tP_n(t)=A_{n+1}P_{n+1}(t)+B_nP_n(t)+A_n^*P_{n-1}(t),
\qquad n\ge 0,\label{rekm}\ee
 where $A_n$ are nonsingular matrices, and $B_n$ are
Hermitian. Conversely, the analogue of Favard's theorem is
also true: if a sequence of matrix polynomials
$\{P_n\}$ of corresponding degree $n=0,1,2,\ldots$, satisfy
(\ref{rekm}) with nonsingular $A_n$ and Hermitian $B_n$,
then there is a positive definite measure matrix $W$ such
that the $P_n$ are orthonormal with respect to $W$.

The three-term recurrence formula easily yields the
 Christoffel-Darboux formula:
\[(w-z)\sum_{k=0}^{n-1}P_k^*(z)P_k(w)=P_{n-1}^*(z)A_nP_n(w)-P_n^*(z)A_n^*P_{n-1}(w),\]
from which for example  it follows that
\[P_{n-1}^*(z)A_nP_n(z)-P_n^*(z)A_n^*P_{n-1}(z)=0,\]
\[\sum_{k=0}^{n-1}P_k^*(z)P_k(z)=P_{n-1}^*(z)A_nP_n'(z)-P_n^*(z)A_n^*P_{n-1}'(z).\]

The orthogonal polynomials $Q_n$ of the second kind
\[Q_n(t)=\int \frac{P_n(t)-P_n(x)}{t-x}dW(x),\qquad n=1,2,\ldots,\]
also satisfy the same recurrence and are orthogonal with respect to
some other matrix measure. For them we have
\[P_{n-1}^*(t) A_n Q_n(t)-P_n^*(z)A_n^* Q_{n-1}(t)\equiv I,\]
and
\[Q_n(t)P_{n-1}^*(t)-P_n(t)Q_{n-1}^*(t)\equiv A_n^{-1}.\]

With the recurrence coefficient matrices $A_n$, $B_n$ one can
form the \dword{block Jacobi matrix}
\[J= \left(\begin{array}{ccccc}
  B_0 & A_0 & 0 & 0 &\cdots  \\
  A_0^* & B_1 & A_1 & 0 & \cdots \\
  0 & A_1^* & B_2 & A_2 & \cdots \\
  0 & 0 & A_2^* & B_2 &\cdots \\
  \vdots&\vdots&\vdots&\vdots&\ddots
\end{array}\right).\]
The zeros of $P_n$ are real and they are
the eigenvalues (with the same multiplicity)
of the $N$-truncated block Jacobi matrix (which is of size $nN$). If
$a$ is a zero then its multiplicity $p$ is at most $N$, the rank of
$P_n(a)$ is $N-p$, and the space of those vectors $v$ for which $P_n(a)v=0$
is of dimension $p$. If we write $x_{n,k}$, $1\le k\le m$, for the different zeros of
$P_n$, and $l_k$ is the multiplicity of $x_{n,k}$, then the matrices
\[\G_k=\frac{1}{({\rm det}(P_n(t)))^{(l_k)}(x_{n,k})}\left({\rm Adj}(P_n(t))
\right)^{(l_k-1)}(x_{n,k})Q_n(x_{n,k}),\quad 1\le k\le m\]
are positive semidefinite of rank $l_k$, and with them
the \dword{matrix quadrature formula}
\[\int P(t)dW(t)=\sum_{k=1}^m P(x_{n,k})\G_{n,k}\]
holds for every matrix polynomial $P$ of degree at most $2n-1$.

If we assume that $A_n\to A$, $B_n\to B$ where $A$ is non-singular,
then
\[P_n(z)P_{n-1}^{-1}(z)A_n^{-1}\to \int\frac{dW_{A,B}(t)}{z-t}\]
locally uniformly outside the cluster set of the zeros, where $W_{A,B}(t)$
is the measure matrix of orthogonality for the sequence of matrix orthogonal 
polynomials $S_n$ with recurrence coefficients $A,B$ for all $n$, i.e., which satisfy the
three-term recurrence
\[ tS_n(t)=A^*S_{n+1}(t)+BS_n(t)+AS_{n-1}(t).
\]
The distribution of the zeros themselves will be $1/N$-times the 
trace of the matrix measure of orthogonality for
another sequence of matrix orthogonal polynomials $R_n$ satisfying
\[ tR_n(t)=AR_{n+1}(t)+BR_n(t)+AR_{n-1}(t),\qquad n\ge 2 ,\]
with appropriate modifications for $n<2$.

\subsubsection{Families of orthogonal polynomials}
If the matrix of orthogonality is diagonal (or similar to a diagonal matrix) with diagonal
entries $\mu_i$,
then the orthogonal matrix polynomials are also diagonal with $i$-th entry equal to
$p_n(\mu_i)$, the $n$-th orthogonal polynomial with respect to $\mu_i$. Many matrix
orthogonal polynomials in the literature can be reduced to
this scalar case. Recently however, some remarkably rich non-reducible families have
been obtained by A. Duran and F. Gr\"unbaum (see \cite{duran1} and the references therein),
which may play the role of the classical orthogonal polynomials in higher dimension. They found
families of matrix orthogonal polynomials that satisfy second order (matrix) differential
equations just like the classical orthogonal polynomials. Their starting point was
a symmetry property between the orthogonality measure matrix and a second order
differential operator. They worked out
several explicit examples. Here is
one of them: $N=2$, the measure matrix (more precisely its
density) is
\[H(t):=e^{-t^2}\left(\begin{array}{cc}
  1+|a|^2 t^4 & at^2\\
  \overline a t^2 & 1
\end{array}\right),\qquad \ t\in \R,\]
where $a\in \C\setminus \{0\}$ is a free parameter.
The corresponding $P_n(t)$
satisfies
\bean P_N''(t)+P_n'(t)\left(\begin{array}{cc}
  -2t & 4at \\
  0 & -2t
\end{array}\right)
&+&P_n(t)\left(\begin{array}{cc}
  -4 & 2a \\
  0 & 0
\end{array}\right)\\[12pt]
&=&\left(\begin{array}{cc}
  -2n-4 & 2a(2n+1) \\
  0 & -2n
\end{array}\right)P_n(t).\eean
There is an explicit Rodrigues' type representation for the polynomials themselves,
and the three-term recurrence (\ref{rekm}) holds with $B_n=0$,
\[A_{n+1}:=\sqrt{\frac{n+1}{2}}\left(\begin{array}{cc}
  \g_{n+3}/\g_{n+2} & a\g_{n+2}\g_{n+1} \\
  0 & \g_n/\g_{n+1}
\end{array}\right),\]
where
\[\g_n^2:=1+\frac{|a|^2}2 {n\choose 2}.\]

\medskip

\subsubsection{Connection with higher order
scalar recurrence}
Matrix orthogonality is closely connected to  \dword{$(2N+1)$-term recurrence}s
for scalar polynomials. To describe this we need the following operators
on polynomials $p$: if $p(t)=\sum_k a_kt^k$, then
\[R_{N,m}(p)=\sum_s a_{sN+m}t^s,\]
i.e., from a polynomial the operator $R_{N,m}$ takes those powers where
the exponent is congruent to $m$ modulo $N$, removes the common factor $t^m$ and
changes $t^N$ to $t$.

Now suppose that $\{p_n\}_{n=0}^\i$ is a sequence of scalar polynomials
of corresponding degree $n=0,1,\ldots$, and suppose that this sequence
satisfies a $(2N+1)$-term recurrence relation
\[t^Np_n(t)=c_{n,0}p_n(t)+\sum_{k=1}^N(\overline{c_{n,k}}p_{n-k}(t)+c_{n+k,k}p_{n+k}(t)),\]
where $c_{n,0}$ is real, $c_{n,N}\not=0$ (and $p_k(t)\equiv 0$ for $k<0$). Then
\[P_n(t)=\left(\begin{array}{ccc}
  R_{N,0}(p_{nN}) & \cdots & R_{N,N-1}(p_{nN}) \\
   R_{N,0}(p_{nN+1}) & \cdots &   R_{N,N-1}(p_{nN}+1)\\
  \vdots & \ddots & \vdots \\
R_{N,0}(p_{nN+N-1}) & \cdots &   R_{N,N-1}(p_{nN}+N-1)
\end{array}\right)\]
is a sequence of matrix orthogonal polynomials with respect
to a positive definite measure matrix. Conversely, if
$P_n=(P_{n,m,j})_{m,j=0}^{N-1}$ is a sequence of orthonormal
matrix polynomials, then the scalar polynomials
\[p_{nN+m}(t)=\sum_{j=0}^{N-1} t^j P_{n,m,j}(t^N),
\qquad 0\le m<N, \ n=0,1,2,\ldots,\]
satisfy a $(2N+1)$-recurrence relation of the above form.

Bolyai Institute,
University of Szeged

Szeged

Aradi v. tere 1, 6720

Hungary, and
\medskip

Department of Mathematics

University of South Florida

Tampa, FL, 33620

USA

\endddoc
\begin{thebibliography}{99}

\bibitem{akhiezer1}  N. I. Ahiezer,
Orthogonal polynomials on several intervals,
{\it   Dokl. Akad. Nauk SSSR}, {\bf  134}(1960), 9--12. (Russian); 
translated as {\it Soviet Math. Dokl.}, {\bf   1}(1960),  989--992. 

\bibitem{akhiezer} N. I. Akhiezer,
{\it The Classical Moment Problem and Some Related Questions in Analysis},
Oliver and Boyd, Edinburgh, 1965 


\bibitem{alsalam} W. A. Al-Salam, Characterization theorems for orthogonal
polynomials, in: Orthogonal Polynomials: Theory and Practice,
M. Ismail and P. Nevai (eds), NATO ASI Series C, {\bf 294}, Kluwer, Dordrecht, 1990. 

\bibitem{alfaro1} M. P. Alfaro, M. Bello, J. M. Montaner and J. L. 
Varona, Some asymptotic properties for orthogonal polynomials with
respect to varying measures, {\it J. Approx. Theory}, {\bf
135}(2005), 22--34. 

\bibitem{alfaro} M. P. Alfaro and L. Vigil, Solution of a problem
of P. Tur\'an, {\it J. Approx. Theory}, {\bf 53}(1988), 195--197. 

\bibitem{Askey} G. E. Andrews, R. Askey and R. Roy, {\it Special Functions},
Encyclopedia of Mathematics and its Applications, Cambridge
University Press, New York, 1999. 

\bibitem{aptekarev1} A. I. Aptekarev,
Asymptotic properties of polynomials orthogonal on a system of contours, and
periodic motions of Toda chains. (Russian)
{\it Mat. Sb. (N.S.)}, {\bf 125(167)}(1984),  231--258. 

\bibitem{aptekarev} A. I. Aptekarev,
Multiple orthogonal polynomials, {\it J. Comput. Appl. Math.}, {\bf
99}(1988), 423--447. 

\bibitem{Aptekarevnikishin} A. I. Aptekarev and E. M. Nikishin,
The scattering problem for a discrete Sturm-Liouville operator,
{\it Mat. Sb.} (N.S.),  {\bf 121}(1983),
 327--358. (Russian)

\bibitem{aptekarev2} A. I. Aptekarev and W. Van Assche, Scalar and matrix
Riemann--Hilbert approach to the strong asymptotics of Pad\'e
approximants and complex orthogonal polynomials with varying
weights, {\it J. Approx. Theory}, {\bf 129}(2004), 129--166. 

\bibitem{Askey1} R. Askey, A transplantation theorem for Jacobi series.  
{\it Illinois J. Math.}, {\bf 13}(1969), 583--590.

\bibitem{Baratchart} L. Baratchart,  R. K\"ustner and V. Totik,
Zero distribution via orthogonality, {\it Annales de l'Institut
Fourier},  {\bf 55}(2005), 1455-1499.

\bibitem{barrios} D. Barrios and G. L\'opez, Ratio asymptotics for
polynomials orthogonal on arcs of the unit circle, {\it Constr. 
Approx.}, {\bf 15}(1999), 1--31. 

\bibitem{bello} M. Bello and G. L\'opez, Ratio and relative
asymptotics of polynomials orthogonal on an arc of the unit
circle, {\it J. Approx. Theory}, {\bf 92}(1998), 216--244. 

\bibitem{blumenthal} O. Blumenthal, {\it \"Uber die Entwicklung einer
willk\"urlichen Funktion nach den Nennern des Kettenbruches f\"ur
$\int_{-\infty}^0 [\phi(\xi)/(z-\xi)]\, d\xi$,} Inaugural Dissertation,
G\"ottingen, 1898. 


\bibitem{bochner} S. Bochner, \"Uber Sturm-Liouvillesche Polynomsysteme,
{\it Math. Zeitschrift}, {\bf 29}(1929),
730--736. 

\bibitem{ismail5} J. A. Charris and M. E. H. Ismail, Sieved orthogonal polynomials. VII. Generalized polynomial mappings,
{\it   Trans. Amer. Math. Soc.},  {\bf 340}(1993), 71--93. 

\bibitem{chihara} T. S. Chihara, {\it An introduction to orthogonal polynomials.}
 Mathematics and its Applications, Vol. {\bf 13},
  Gordon and Breach Science Publishers, New York-London-Paris, 1978. 

\bibitem{damems} E. Daems and A. B. J. Kuijlaars, A Christoffel-Darboux formula
for multiple orthogonal polynomials, {\it J. Approx. Theory}, {\bf 130}(2004), 190--202. 

\bibitem{Deift1} P. Deift, {\it Orthogonal Polynomials and Random Matrices:
A Riemann--Hilbert Approach}, Courant Lecture Notes, 1999. 


\bibitem{simondely} F. Delyon, B. Simon, and B. Souillard, From power pure point to
continuous spectrum in disordered systems, {\it Annales Inst. 
Henri Poincar\'e}, {\bf 42}(1985), 283--309. 

\bibitem{denisov} S. Denisov, On Rakhmanov's theorem for Jacobi matrices,
{\it Proc. Amer. Math. Soc.}, {\bf 132}(2004), 847--852. 

\bibitem{denissimon} S. Denisov and B. Simon, Zeros of orthogonal polynomials on the real line,
{\it J. Approx. Theory},  {\bf 121}(2003), 357--364. 

\bibitem{xu} C. F. Dunkl and Y. Xu, {\it Orthogonal polynomials of several variables},
Encyclopedia
of Mathematics and its Applciations, Cambridge University Press, New York, 2001. 

\bibitem{duran1} A. J. Dur\'an and F. A. Gr\"unbaum, A survey on orthogonal
matrix polynomials satisfying second order differential equations,
{\it J. Comput. Appl. Math.}, {\bf 178}(2005), 169-190. 

\bibitem{Its} A. S. Fokas, A. R. Its and A. V. Kitaev, The isomonodromy
approach to matrix models in 2D quantum gravity,
{\it Comm. Math. Phys.}, {\bf 147}(1992), 395--430. 

\bibitem{freud} G. Freud, {\it Orthogonal Polynomials.}
Akad\'emiai
Kiad\'o/Pergamon Press, Budapest, 1971. 


\bibitem{gautschi} W. Gautschi, {\it Orthogonal Polynomials:
Computation and Approximation}, Numerical Mathematics and
Scientific Computation, Oxford University Press, New York, 2004. 

\bibitem{Kuijl} W. Gautschi and A. B. J. Kuijlaars, Zeros and critical points
of Sobolev orthogonal polynomials, {\it J. Approx. 
Theory}, {\bf 91}(1997), 117--137. 

\bibitem{gerlubmarc} J. Geronimo, D. S. Lubinsky and F. Marcellan, Asymptotics for
Sobolev orthogonal polynomials for exponential weights, to appear in
{\it Constructive Approx.}

\bibitem{geronimo} J. S. Geronimo and W. Van Assche, Orthogonal polynomials on
several intervals via a polynomial mapping, {\it Trans. Amer. Math. Soc.},
{\bf 308}(1988), 559--581. 

\bibitem{geronimus} Ya. L. Geronimus, On the character of the
solution of the moment problem in the case of the periodic in the
limit associated fraction, {\it Bull. Acad. Sci. USSR Math.}, {\bf
5}(1941), 203--210. [Russian]

\bibitem{geronimus1} Ya. L. Geronimus, On polynomials orthogonal
on the circle, on trigonometric moment problem, and on allied
Carath\'eodory and Schur functions, {\it Mat. Sb.}, {\bf 15}(1944),
99--130. [Russian]

\bibitem{geronimus2} Ya. L. Geronimus, On the trigonometric moment
problem, {\it Ann. of Math.}, {\bf 47}(1946), 742--761. 

\bibitem{geronimus3} Ya. L. Geronimus, {\it Orthogonal polynomials on the 
unit circle and their applications}, Amer. Math. Soc. Translation, 1954,
Vol. 104. 


\bibitem{goltotik} L. Golinskii and V. Totik, Orthogonal polynomials
from Stieltjes to Simon, Proc. Conf. honoring B. Simon,
AMS Conference Series (to appear).


\bibitem{Goncharrahmanov} A. A. Gonchar and E. A. Rakhmanov,
 On the convergence of simultaneous Pad\'e approximants for systems
 of functions of Markov type, {\it Number theory, mathematical
 analysis and their applications. Trudy Mat. Inst. Steklov.},
 {\bf   157}(1981), 31--48. (Russian)

\bibitem{Goncharsuetin} A. A. Gonchar and S. P. Suetin, On
Pad\'e  approximation of Markov type meromorphic functions,
{\it Current Problems in Mathematics,
5.}, Rossiiskaya Akademiya Nauk, Matematicheskii Institut im. V. A. Steklova,
Moscow, 2004. 68 pp. (electronic). (Russian)


\bibitem{Ismail} M. E. H. Ismail, {\it Classical and quantum orthogonal polynomials
in one variable}, Encyclopedia in Mathematics, Cambridge University
Press, 2005. 

\bibitem{kalyagin} V. Kaliaguine, A note on the asymptotics of orthogonal
polynomials on a complex arc: the case of a measure with a discrete part. ,
{\it J. Approx. Theory}, {\bf  80}(1995),  138--145. 


\bibitem{khrush3} S. Khrushchev, A singular Riesz product in the
Nevai class and inner functions with the Schur parameters in
$\cap_{p>2} \ell^p$, {\it J. Approx. Theory}, {\bf 108}(2001),
249--255. 

\bibitem{Kruschev} S. Khrushchev, Continued fractions and orthogonal polynomials
on the unit circle, Proc. Conf. OPSFA, Copenhagen, 2003. 

\bibitem{killip} R. Killip and B. Simon,
 Sum rules for Jacobi matrices and their applications to spectral theory,
 {\it Annals of Math.}, {\bf 158}(2003), 253-321

\bibitem{Kriecher} T. Kriecherbauer and K. T-R. McLaughlin, Strong asymptotics
of polynomials orthogonal with respect to a Freud weight, {\it
Internat. Math. Research Notices}, {\bf 6}(1999), 299--333. 

\bibitem{Levlub1} A. Levin and D. S. Lubinsky,
{\it Christoffel Functions and Orthogonal Polynomials for Exponential
Weights on $[-1,1]$}, Memoirs Amer. Math. Soc. , {\bf 535}, (1994)

\bibitem{Levlub} A. Levin and D. S. Lubinsky,
{\it Orthogonal polynomials for exponential weights},
 CMS Books in Mathematics/Ouvrages de Mathïmatiques de la SMC, {\bf 4},
  Springer-Verlag, New York, 2001. 

\bibitem{Lopez} G. L\'opez, H. Pijeira-Cabrera and I. P. Izquierdo,
Sobolev orthogonal polynomials in the complex plane,
{\it J. Comput. Appl. Math.}, {\bf 127}(2001), 219-230. 

\bibitem{duran} P. L\'opez-Rodriguez and A. J. Dur\'an, Orthogonal matrix polynomials,
Laredo Lecture Notes (Lecture notes of the 1st SIAG Summer School on
Orthogonal Polynomials and Special Functions, Laredo, Spain, July 24-29
2000). Nova Science Publisher (2003). R. \'Alvarez-Nodarse, F. Marcell\'an
and W. Van Assche (Eds. )

\bibitem{Lub} D. S. Lubinsky, {\it Strong Asymptotics for Extremal Errors
Associated with Erd\"os-type Weights}, Pitman Res. Notes, Vol. {\bf 202}, Longmans,
Nurn Mill, Essex, 1989. 

\bibitem{Lubsaff} D. S. Lubinsky and E. B. Saff,
{\it Strong Asymptotics for Extremal Polynomials Associated with Weights on $(-\i,\i)$},
Springer Lecture Notes, Vol. {\bf 1305}, Springer Verlag, Berlin, 1988. 

\bibitem{Marcellan} F. Marcell\'an, M. Alfaro and M. L. Rezola, Orthogonal polynomials
on Sobolev spaces: old and new directions, {\it J. Comput. Appl. 
Math.,} {\bf 48}(1993), 113--131. 

\bibitem{Martinez1} A. Martinez-Finkelshtein, Analytic aspects of Sobolev orthogonal
polynomials, {\it J. Comput. Appl. Math.}, {\bf 99}(1998),
491--510. 

\bibitem{Martinez} A. Martinez-Finkelshtein, Analytic properties of Sobolev orthogonal
polynomials revisited, {\it J. Comput. Appl. Math.}, {\bf
127}(2001), 255-266. 


\bibitem{martinez1} A. Martinez-Finkelshtein, \SZ\ polynomials: a
view from the Riemann--Hilbert window, preprint
arXiv:math. CA/0508117, 2005. 

\bibitem{mastroianni} G. Mastroianni and G. Milovanovic,
{\it Approximation of Functions and Functionals with Application},
World Scientific,  to appear. 

\bibitem{mastto} G. Mastroianni and V. Totik,  Polynomial inequ\-alities
with do\-u\-b\-ling and $A_\infty$ weights, {\it Con\-structive
Ap\-p\-rox.}, {\bf 16}(2000), 37--71. 

\bibitem{mnt1} A. M\'at\'e, P. Nevai and V. Totik, Extensions of
Szeg\H o's theory of orthogonal polynomials, {\it
Constructive Approx.}, {\bf 3}(1987), 51--72. 

\bibitem{mnt2} A. M\'at\'e,  P. Nevai and V. Totik,  Strong and weak
convergence of orthogonal polynomials, {\it Amer. J. Math.},
{\bf 109}(1987), 239--282. 

\bibitem{mnt3} A. M\'at\'e,  P. Nevai and V. Totik, Szeg\H o's
extremum problem on the unit circle, {\it Annals of Math.},
{\bf 134}(1991), 433--453. 

\bibitem{mehta} M. Mehta, {\it Random matrices}, second ed.,
Academic Press, Inc. , Boston, 1991. 

\bibitem{mhaskar} H. N. Mhaskar and E. B. Saff, On the distribution
of zeros of polynomials orthogonal on the unit circle, {\it J. 
Approx. Theory}, {\bf 63}(1990), 30--38. 

\bibitem{Nevai} P. Nevai, {\it Orthogonal Polynomials,}
Memoirs  Amer. Math. Soc., {\bf 213}(1979)

\bibitem{nevtot1} P. Nevai and V. Totik, Orthogonal polynomials and
their zeros, {\it Acta Sci. Math.}(Szeged), {\bf 53}(1989),
99--104. 

\bibitem{nevtot2} P. Nevai and V. Totik, Denisov's theorem on
recurrent coefficients, {\it J. Approx. Theory}, {\bf 127}(2004),
240--245. 

\bibitem{nikishin} E. M. Nikishin, An estimate for orthogonal
polynomials, {\it Acta Sci. Math.}(Szeged), {\bf 48}(1985),
395--399. [Russian]


\bibitem{Nuttall1} J. Nuttall, Asymptotics for diagonal
Hermite-Pad\'e approximants, {\it J. Approx. Theory},
{\bf 42}(1984), 299--386. 

\bibitem{Nuttall2} J. Nuttall, Pad\'e polynomial  asymptotics from
a singular integral equation, {\it Constructive Approx.}, {\bf 6}(1990), 157--166. 

\bibitem{Pastur} L. Pastur and A. Figotin, {\it Spectra of Random and Almost-Periodic
 Operators}, Grundlehren der mathematischen Wissenschaften, Vol. {\bf 297}, Springer Verlag,
Berlin, 1992. 

\bibitem{peher2} F. Peherstorfer,  On Bernstein-Szeg\H o orthogonal polynomials
on several intervals,
{\it SIAM J. Math. Anal.}, {\bf 21}(1990), 461--482. 

\bibitem{peher3} F. Peherstorfer,  On Bernstein-Szeg\H o orthogonal polynomials
on several intervals. II. 
Orthogonal polynomials with periodic recurrence coefficients,
{\it   J. Approx. Theory},  {\bf 64}(1991),   123--161. 

\bibitem{peher5} F. Peherstorfer,  Orthogonal and extremal polynomials on several intervals,
{\it  J. Comput. Appl. Math.},  {\bf 48}(1993),  187--205. 

\bibitem{peher4} F. Peherstorfer, Zeros of polynomials orthogonal on several intervals. 
{\it Int. Math. Res. Notes},  {\bf 7}(2003),  361--385. 

\bibitem{peher1} F. Peherstorfer and R. Steinbauer, On polynomials orthogonal
on several intervals,
Special functions (Torino, 1993),  {\it Ann. Numer. Math.}, {\bf
2}(1995), 353--370. 

\bibitem{rahmanov1} E. A. Rakhmanov, On V. A. Steklov's problem in the
theory of orthogonal polynomials, (Russian)
{\it Dokl. Akad. Nauk SSSR},  {\bf 254}(1980),  802--806. 

\bibitem{rakhmanov} E. A. Rakhmanov, On the asymptotics of the ratio of
orthogonal polynomials, {\it Math. USSR Sbornik,} {\bf 46}(1983), 105--117;
Russian original, {\it Mat. Sb.}, {\bf 118}(1982), 104--117. 

\bibitem{safftotik2} E. B. Saff and V. Totik, 
What parts of a measure's support attract zeros of the corresponding orthogonal polynomials?
{\it Proc. Amer. Math. Soc.},  {\bf 114}(1992), 185--190. 

\bibitem{safftot} E. B. Saff and V. Totik,
{\it Logarithmic Potentials with External Fields}, 
Grund\-leh\-ren der mathematischen {W}issenschaften,
{\bf 316}, Springer-Verlag, Berlin, 1997. 

\bibitem{Simonmoment} B. Simon, The classical moment problem as a self-adjoint finite
difference operator, {\it Advances in Math.}, {\bf 137}(1998),
82--203. 

\bibitem{simon6} B. Simon, Ratio asymptotics and weak asymptotic
measures for orthogonal polynomials on the real line, {\it J. 
Approx. Theory}, {\bf 126}(2004), 198--217. 

\bibitem{simon1} B. Simon, {\it Orthogonal Polynomials on the Unit
Circle}, V. 1: Classical Theory, AMS Colloquium Series, American
Mathematical Society, Providence, RI, 2005. 

\bibitem{simon2} B. Simon, {\it Orthogonal Polynomials on the Unit
Circle}, V. 2: Spectral Theory, AMS Colloquium Series, American
Mathematical Society, Providence, RI, 2005. 


\bibitem{simzlat} B. Simon and A. Zlato\v{s}, Sum rules and the \SZ\
condition for orthogonal polynomials on the real line, {\it Comm. 
Math. Phys.}, {\bf 242}(2003), 393--423. 

\bibitem{zlatos} B. Simon and A. Zlato\v{s}, Higher-order \SZ\
theorems with two singular points, {\it J. Approx. Theory}, {\bf
134}(2005), 114--129. 

\bibitem{stahl3} H. Stahl, Divergence of diagonal Pad\'e approximants and the
asymptotic behavior of orthogonal polynomials associated with nonpositive measures,
{\it Constr. Approx.}, {\bf  1}(1985),  249--270. 

\bibitem{stahl2} H. Stahl,  Orthogonal polynomials with complex-valued weight function. 
I, II. {\it Constr. Approx.}, {\bf  2}(1986), 225--240, 241--251. 

\bibitem{stahl1} H. Stahl, Orthogonal polynomials with respect to complex-valued measures. 
Orthogonal polynomials and their applications (Erice, 1990),
{\it  IMACS Ann. Comput. Appl. Math.}, {\bf 9}(1991),  139--154. 

\bibitem{stahltotik} H. Stahl and V. Totik, {\it General Orthogonal
Polynomials}, Encyclopedia of Mathematics and its Applications, {\bf 43},
Cambridge University Press, New York, 1992.

\bibitem{Suetin1}   S. P. Suetin, On the uniform convergence of diagonal
Pad\'e approximants for hyperelliptic functions, {\it Mat. Sb.},
{\bf  191}(2000),  81--114; (Russian),
 translation in  {\it Sb. Math.},  {\bf 191}(2000), 1339--1373.

\bibitem{Suetin2} S. P. Suetin,  Pad\'e approximants and the
effective analytic continuation of a power series, {\it  Uspekhi
Mat. Nauk,} {\bf  57}(2002), 45--142; (Russian), translation in
{\it Russian Math. Surveys},  {\bf 57}(2002),  43--141.


\bibitem{szabados} J. Szabados and P. V\'ertesi, {\it Interpolation of Functions},
World Scientific, Singapore, 1990. 

\bibitem{Szego} G. Szeg\H{o}, {\it Orthogonal Polynomials},
Coll. Publ. , {\bf XXIII}, Amer. Math. Soc. , Providence, 1975. 

\bibitem{totik1} V. Totik, Orthogonal polynomials with ratio asymptotics, {\it
Proc. Amer. Math. Soc.}, {\bf 114}(1992), 491--495. 

\bibitem{totik2} V. Totik, Asymptotics for Christoffel functions for general
measures on the real line, {\it J. D{'}Analyse Math.}, {\bf
81}(2000), 283--303. 

\bibitem{Assche} W. Van Assche, {\it Asymptotics for Orthogonal Polynomials},
Spriger Lecture Notes, Vol. {\bf 1265}, Springer Verlag, Berlin, 1987. 

\bibitem{assche2} W. Van Assche, Hermite-Pad\'e approximation and multiple
orthogonal polynomials, preprint. 

\bibitem{assche1} W. Van Assche, J. S. Geronimo, A. B. J. Kuijlaars,
Riemann--Hilbert problems
for multiple orthogonal polynomials, in: J. Bustoz et al. (eds), Special Functions 2000:
Current Perspectives and Future Directions, Kluwer, Dordrecht, 2001, 23--59. 

\bibitem{widom} H. Widom, Extremal polynomials associated with a system of
curves in the complex plane, {\it Adv. Math.}, {\bf 3}(1969), 127--232. 




\end{thebibliography}
